\documentclass[11pt]{amsart} 
\usepackage{amsmath,amssymb,a4wide} 
\usepackage{color}
\usepackage[latin1]{inputenc}
\usepackage{color}
\usepackage[dvipsnames]{xcolor}
\usepackage{mathrsfs}

\newtheorem{theorem}{Theorem} 
\newtheorem{lem}[theorem]{Lemma} 
\newtheorem{prop}[theorem]{Proposition} 
\newtheorem{rem}[theorem]{Remark} 
\newtheorem{cor}[theorem]{Corollary}
\newtheorem{defi}[theorem]{Definition}

\newcommand{\pr}{{\bf \textit{Proof : }}}

\newcommand{\cqfd}{{\nobreak\hfil\penalty50\hskip2em\hbox{}\nobreak\hfil
$\square$\qquad\parfillskip=0pt\finalhyphendemerits=0\par\medskip}}

\newcommand{\R}{{\mathbb R}} 
 
\newcommand{\TD}{{\mathbb T}}

\let\eps\varepsilon

\DeclareMathOperator*{\dive}{div}

\title[On Navier-Stokes-Korteweg and Euler-Korteweg Systems]{On Navier-Stokes-Korteweg and Euler-Korteweg Systems: Application to Quantum Fluids MOdels}

\author[D. Bresch]{D. Bresch}
\address{Laboratoire de Math\'ematiques, CNRS UMR 5127,
Universit\'e  Savoie Mont-Blanc, 73376 Le Bour\-get-du-Lac, France; 
e-mail: Didier.Bresch@univ-smb.fr}

\author[M. Gisclon]{M. Gisclon}
\address{Laboratoire de Math\'ematiques, CNRS UMR 5127,
Universit\'e Savoie Mont-Blanc, 73376 Le Bour\-get-du-Lac, France; 
e-mail: gisclon@univ-smb.fr}

\author[I. Violet]{I. Lacroix-Violet} 
\address{Laboratoire de Math\'ematiques, CNRS UMR 8524 
Universit\'e de Lille 1, Villeneuve d'Ascq, France;
e-mail: ingrid.violet@univ-lille1.fr}

\begin{document}

\maketitle

\begin{abstract}

 In this paper, the main objective is to generalize to the Navier-Stokes-Korteweg (with 
 density dependent viscosities satisfying the BD relation) and   Euler-Korteweg systems  a recent relative entropy [proposed by D. {\sc Bresch}, P. {\sc Noble} and J.--P.~{\sc Vila}, (2016)] introduced for the compressible Navier-Stokes equations with a  li\-near density dependent  shear viscosity and a zero bulk viscosity.   
  As a concrete application, this helps to justify mathematically the convergence between global weak solutions of the quantum Navier-Stokes system [recently obtained simultaneously by  {\sc I.} {\sc Lacroix-Violet} and  {\sc A.} {\sc Vasseur} (2017)]  and dissipative solutions of the quantum Euler system when the viscosity coefficient tends to zero: This selects a dissipative solution as the limit of a viscous system.  We also get weak-strong uniqueness for the Quantum-Euler and for the Quantum-Navier-Stokes equations. 
  Our results are based on the fact that Euler-Korteweg systems and correspon\-ding Navier--Stokes-Korteweg systems can be reformulated through an augmented  system such as the compressible Navier-Stokes system with density dependent viscosities satisfying the BD algebraic relation. This was also observed recently [by D. {\sc Bresch}, F. {\sc Couderc}, P. {\sc Noble} and J.--P. {\sc Vila}, (2016)] for the Euler-Korteweg system for numerical purposes.
 As a by-product of our analysis, we show that this augmented formulation helps to define relative entropy estimates for the Euler-Korteweg  systems in a simplest way  compared to recent works [See {\sc  D. Donatelli, E. Feireisl, P.~Marcati} (2015) and  {\sc J.~Giesselmann, C. Lattanzio, A.-E. Tzavaras} (2017)] with less hypothesis required 
 on the capillary coefficient.  
   
\end{abstract}

\bigskip

{\small\noindent 

} 


\section{Introduction}
Quantum fluid models have attracted a lot of attention in the last decades due to te variety of applications. Indeed, such models can be used to describe superfluids \cite{LoMo93}, quantum semiconductors \cite{FeZhou93}, weakly interacting Bose gases \cite{Grant73} and quantum trajectories of Bohmian mechanics \cite{Wyatt05}. Recently some dissipative quantum fluid models have been derived. In par\-ticular,  under some assumptions and using a Chapman-Enskog expansion in Wigner equation, the authors have obtained in \cite{BruMe09} the so-called quantum Navier-Stokes model. Roughly spea\-king, it corresponds to the classical Navier-Stokes equations with a quantum correction term.  The main difficulties of such models lie in the highly nonlinear structure of the third order quantum term and the proof of positivity (or non-negativity) of the particle density. Note that formally, the quantum Euler system corresponds to the limit of the quantum Navier-Stokes model when the viscosity coefficient tends to zero. This type of models belong to more general classes of models: the Navier-Stokes-Korteweg and the Euler-Korteweg systems. Readers interested by Korteweg type systems are referred to the following articles and books: \cite{Ko, Va, CaHi, DuSe, RoWi, Ni, HeMa} and references cited therein.

The goal of this paper is to extend to these two Korteweg systems a recent relative entropy proposed in \cite{BNV} introduced for the compressible Navier-Stokes equations with a linear density dependent shear viscosity and a zero bulk viscosity. This leads for each system to the definition of what we call a dissipative solution following the concept introduced by P.-L. Lions in the incompressible setting (see \cite{L}) and later extended to the compressible framework (see \cite{FeNo, FNL, BN, S} for constant viscosities and \cite{Ha, BNV} for density dependent viscosities). As a consequence we obtain some weak-strong uniqueness results and as an application, we can use it to show that a global weak solution (proved in \cite{LaVa}, which is also a dissipative one) of the quantum Navier-Stokes system converges to a dissipative solution of the quantum Euler system.    Our results will be compared to recent results in \cite{DFM, GLT} showing that we relax one hypothesis on the capillarity coefficient by introducing entropy-relative solutions of an augmented system.  Note also the interesting paper \cite{AnSp} where the authors prove the existence of global weak solutions of the quantum-Navier-Stokes equations with a different method compared to \cite{LaVa}. By the way we cannot use such global weak solutions because  capillarity and viscosity  magnitudes are linked together in their study.
Let us also the interesting new paper \cite{CaCaHi} where the authors investigate the long-time behavior of solutions to the isothermal Euler-Korteweg system.
\medskip

Let us now present in more details the models of interest here. Note that for the convenience of the reader all the operators are defined in Section \ref{annex_operators}. Let  $\Omega=\TD^d$ be the torus in dimension $d$ (in this article $1 \leq d \leq 3$). 

\medskip

\noindent {\bf Euler-Korteweg system.} Following the framework of the paper, we first present the Euler-Kortewg system and then the Navier-Stokes Korteweg one.  Note that in all the paper, the systems are supplemented with the following initial conditions
\begin{equation}
\label{condI}
\rho\vert_{t=0}=\rho_0, \,  \quad (\rho \, u)\vert_{t=0}= \rho_0 u_0 \quad \hbox{for} ~x \in \Omega.
\end{equation}
with the regularity $\rho_0 \geq 0, \, \rho_0 \in L^{\gamma}(\Omega),
 \, \rho_0 \, |u_0|^2 \in L^1(\Omega),
  \, \sqrt{K(\rho_0)}\nabla \rho_0 \in L^2(\Omega)$.
The Euler-Korteweg system describe the time evolution, for $t>0$ of the density $\rho=\rho(t,x)$ and the momentum $J=J(t,x)=\rho(t,x)u(t,x)$ (with $u$ the velocity), 
for $x \in \Omega$, of an inviscid fluid.  The equations can be written in the form (\cite{DFM}):
\begin{eqnarray} 
&\partial_t\rho+\dive J=0, & \label{EK1} \\
&\partial_t J +\dive \left( \dfrac{J \otimes J}{\rho} \right)+\nabla(p(\rho))=\eps^2 \, \rho \nabla\left(K(\rho)\Delta \rho+\dfrac{1}{2}K'(\rho)|\nabla \rho|^2\right), & \label{EK2}
\end{eqnarray}
where   
$K : (0,\infty) \rightarrow (0,\infty)$ is a smooth function and $p$ is the pressure function given by $p(\rho)=\rho^\gamma$ for $\gamma>1$. Note that it could be interesting to consider non-monotone pressure laws as in \cite{GiTz} and \cite{GLT}. The coefficient $\eps$ stands for the Planck constant.
In this paper we will consider a function $K(\rho)$ which behaves as $ \rho^s$ with $s\in \R$.
  As mentioned in \cite{DFM}, 
$$\rho \nabla\left(K(\rho) \Delta \rho + \dfrac{1}{2}K'(\rho)|\nabla\rho|^2\right)=\dive(\mathbb K),$$
with
$$\mathbb K =\left(\rho \dive(K(\rho)\nabla\rho)+\dfrac{1}{2}(K(\rho)-\rho K'(\rho))|\nabla\rho|^2\right)\mathbb I_{\R^d}-K(\rho)\nabla \rho \otimes \nabla \rho.$$
Observing that $\mathbb K$ may be written 
\begin{equation}\label{Korteweg}
\mathbb K =\left(\dive(\rho K(\rho)\nabla\rho)-\dfrac{1}{2}(K(\rho)+\rho K'(\rho))|\nabla\rho|^2\right)\mathbb I_{\R^d}-K(\rho)\nabla \rho \otimes \nabla \rho.
\end{equation}
and following the ideas of \cite{BCNV} with 
\begin{equation}\label{muK}
\mu'(\rho)=\sqrt{\rho \, K(\rho)},
\end{equation}
we can define the drift velocity $v$ by $$v= \sqrt{\dfrac{K(\rho)}{\rho}} \nabla \rho= \dfrac{\nabla (\mu(\rho))}{\rho}$$
and show the following generalization of the Bohm identity:
$$ \dive(\mathbb K) = {\rm div}(\mu(\rho)\nabla v) + \frac{1}{2}\nabla(\lambda(\rho) {\rm div} v)$$
with 
$$\lambda(\rho)= 2(\mu'(\rho)\rho - \mu(\rho)).$$

\begin{rem}  Note that the relation between $\lambda$ and $\mu$ is exactly the BD relation
found in  \cite{BrDe} in the Navier-Stokes setting: see the  Navier-Stokes-Korteweg part
below. 
\end{rem}

\noindent We will choose $K(\rho)$ as:
$$K(\rho)= \frac{(s+3)^2}{4} \rho^s \hbox{ with } s\in \R  \qquad  \hbox{ in order to get } \qquad 
 \mu(\rho)= \rho^{(s+3)/2}.$$
  This multiplicative constant in the definition of $K$ does not affect any generality, it suffices to change the definition of $\varepsilon$. 
 Then, we  obtain  the following augmented formulation for the  Euler-Korteweg Equations \eqref{EK1}-\eqref{EK2}:
\begin{eqnarray} 
 &\partial_t\rho+\dive (\rho \, u)=0, & \label{EKA1old} \\
 &\partial_t(\rho \, u)+\dive (\rho \, u \otimes u)+\nabla p(\rho)=\eps \,  \left[ \dive ( \mu(\rho)\nabla \bar{v}) + \dfrac{1}{2} \nabla (\lambda(\rho)\dive \bar{v} )\right], &\label{EKA2old} \\ 
&\partial_t(\rho  \, \bar{v})+\dive (\rho \, \bar{v} \otimes u) =  \varepsilon \left[-  \dive ( \mu(\rho) \,  ^t  \nabla u) - \dfrac{1}{2} \nabla( \lambda(\rho) \dive u)\right], \label{EKA3old}
\end{eqnarray}
with 
\begin{equation}\label{EKA4old}
 \lambda(\rho)=2(\rho \, \mu'(\rho)-\mu(\rho)), 
\qquad \bar{v} =  \varepsilon \nabla \mu(\rho)/\rho.
\end{equation}
System \eqref{EKA1old}-\eqref{EKA4old} is called the Euler-Korteweg augmented system in all the sequel. It has been firstly introduced in this conservative form in \cite{BCNV} to propose a useful construction of a numerical scheme with entropy stability property under a  hyperbolic CFL condition for such dispersive PDEs. 
augmented system, the second order operator matrix is skew-symetric.

\medskip

\noindent {\it The Quantum Euler Equations.}
Note that the choice $K(\rho)=1/\rho$ (which gives $\mu(\rho)=\rho$ and $\lambda(\rho)=0$) leads to the Bohm identity
$$\rho \nabla(K(\rho)\Delta \rho+\dfrac{1}{2}K'(\rho)|\nabla \rho|^2)=\dive(\rho \nabla v)=2 \rho \nabla \left( \dfrac{\Delta \sqrt{\rho}}{\sqrt{\rho}} \right).$$
In that case the system \eqref{EKA1old}-\eqref{EKA4old} becomes
\begin{eqnarray} 
 &\partial_t\rho+\dive (\rho \, u)=0, & \label{CQEK1}\\
 &\partial_t(\rho \, u)+\dive (\rho \, u \otimes u)+\nabla(p(\rho))=\eps \, \dive (\rho \,   \nabla \bar{v}), &\label{CQEK2} \\ 
&\partial_t(\rho  \, \bar{v})+\dive (\rho \, \bar{v} \otimes u) = -   \varepsilon \dive (\rho \,  ^t  \nabla u), \label{CQEK3}
\end{eqnarray}
with 
\begin{equation}\label{CQEK4}
\bar{v}= \varepsilon \nabla\log\rho
\end{equation}
which corresponds to the augmented formulation of the quantum Euler system:
\begin{eqnarray}
&\partial_t\rho+\dive (\rho u)=0, & \label{EQ1}\\
&\partial_t(\rho u)+\dive (\rho u \otimes u)+\nabla p(\rho)=2 \, \varepsilon^2 \, \rho \, \nabla\left(\dfrac{\Delta\sqrt{\rho}}{\sqrt{\rho}}\right).
&\label{EQ2}
\end{eqnarray}
Then such a choice gives rise the so called quantum fluid system for which the global existence of weak solutions of \eqref{EQ1}--\eqref{EQ2}  has been shown in \cite{AM09, AM12} and more recently in \cite{CaDaSa} assuming the initial velocity irrotational namely ${\rm curl}(\rho_0 u_0) =0.$ Note that the quantum term is written as \eqref{Korteweg} 
in these papers, namely 
\begin{eqnarray} \label{QQ}
2 \rho \, \nabla\left(\dfrac{\Delta\sqrt{\rho}}{\sqrt{\rho}}\right)
& = & {\rm div}\bigl( \nabla (\rho\nabla \log \rho) 
     -\rho \nabla \log\rho  \otimes \nabla \log \rho \bigr) 
\end{eqnarray}
observing that 
$$\sqrt \rho \nabla \log \rho =  2 \nabla \sqrt \rho.$$
   The existence of local strong solutions has also been proved (see \cite{BDD2007}) and
 global well-posedness for small irrotational data has been performed recently in \cite{AuHa}
 assuming a natural stability condition on the pressure.  We refer to \eqref{CQEK1}-\eqref{CQEK3} as the quantum Euler augmented system in all the paper.
 
 \medskip
\noindent {\it Important remark.} Differentiating in space the mass equation in ${\mathcal D}'((0,T)\times \Omega)$  we get
$$\partial_t \nabla \rho + \nabla \dive (\rho u) = 
     \partial_t \nabla \rho + \dive ({}^t\nabla (\rho u))=0$$
which may be written
$$\partial_t \nabla \rho + \dive (\rho \nabla \log\rho \otimes u)
                                         + \dive\bigl({}^t\nabla(\rho u) - \rho \nabla \log\rho \otimes u)=0$$ 
This formula will be used to show that global weak solutions of the Quantum-Euler
system \eqref{EQ1}--\eqref{EQ2} with the quantum term written as \eqref{QQ} 
 will be global weak solutions of the Quantum-Euler system in its  augmented form.

\medskip

Note that the quantum correction $(\Delta \sqrt{\rho})/\sqrt{\rho}$ can be interpreted as a quantum potential, the so-called Bohm potential, which is well known in quantum mechanics.  This  Bohm potential  arises from the fluid dynamical
formulation of the single-state Schr\"odinger equation. 
The non-locality of quantum mechanics is approximated by the fact that the equations of state do not only depend on the particle density but also on its gradient. 
These equations were employed to model field emissions from metals and steady-state tunneling in metal-insulator-metal structures and to simulate ultra-small semiconductor devices.

\medskip

\noindent {\bf  Navier-Stokes-Korteweg system.}  
Let us consider the compressible Navier-Stokes-Korteweg system with density dependent
viscosities $\mu(\rho)$ and $\lambda(\rho)$ satisfying the  BD relation
$$\lambda(\rho) = 2(\mu'(\rho)\rho - \mu(\rho)),$$
and with the capillarity coefficient $K(\rho)$  linked to the shear viscosity $\mu(\rho)$ in the following manner
$$ K(\rho) = [\mu'(\rho)]^2/\rho \hbox{ with } \mu(\rho)= \rho^{(s+3)/2} \hbox{ with } s\in \R.$$

\begin{rem} With this choice of shear viscosity, the relation between the capillarity
coefficient and the viscosity gives a capillarity coefficient proportional to $\rho^s$.
\end{rem}

Then using the identity given in the Euler-Korteweg part, the Navier-Stokes-Korteweg system can be written for $x \in \Omega$ and  $t>0$,
\begin{eqnarray}
&\partial_t\rho+\dive (\rho u)=0, & \label{NSI0} \\
&\partial_t(\rho u)+\dive (\rho u \otimes u)+\nabla p(\rho)-2 \, \nu  \, \dive (\mu(\rho) D(u)) -  \nu \nabla (\lambda(\rho)\dive u)  \nonumber  \\
 & =
\eps^2\left[ ( \dive(\mu(\rho)^t \nabla v )+\dfrac{1}{2}  \nabla(\lambda(\rho) \dive v)\right],
& \label{NSI1}
\end{eqnarray}
in which the symmetric part of the velocity gradient is $D(u)=\dfrac{1}{2}(\nabla u +^t  \nabla u)$. The parameter $\nu>0$ stands for the viscosity constant. Multiplying \eqref{NSI0} by $\mu'(\rho)$ and taking the gradient, we have the following equation on $v$:
\begin{eqnarray} \label{NSM3}
\partial_t(\rho  \, v)+\dive (\rho \, v \otimes u) +  \dive (\mu(\rho ) \,  ^t  \nabla u) + \dfrac{1}{2} \nabla( \lambda(\rho) \dive u)=0.
\end{eqnarray}
Moreover defining the intermediate velocity, called effective velocity, $w=u+\nu \, v$, equations \eqref{NSI1}  and \eqref{NSM3} lead to
\begin{eqnarray*}
\partial_t(\rho \, w)+\dive (\rho \, w \otimes u)+\nabla(p(\rho))-\nu \,  \dive (\mu(\rho )  \nabla w) -\dfrac{\nu}{2} \nabla( \lambda(\rho) \dive w)   \\
 =
 (\varepsilon^2 -\nu^2) \, [ \dive (\mu(\rho )  \nabla v) + \dfrac{1}{2} \nabla(\lambda(\rho) \dive v)].
\end{eqnarray*}
Then \eqref{NSI0}-\eqref{NSI1} may be reformulated through the following augmented system:
\begin{eqnarray}
&\partial_t\rho+\dive (\rho u)=0, & \label{NSArhoold} \\
&\partial_t(\rho \, w)+\dive (\rho \, w \otimes u)+\nabla(p(\rho))-\nu \,  \dive (\mu(\rho ) \nabla w)  -\dfrac{\nu}{2} \nabla( \lambda(\rho) \dive w) \nonumber   \\
& =
 (\varepsilon^2 -\nu^2)\, [ \dive (\mu(\rho )   \nabla v) + \dfrac{1}{2} \nabla(\lambda(\rho) \dive v)], & \label{NSAwold} \\
&\partial_t(\rho  \, v)+\dive (\rho \, v \otimes u) +  \dive (\mu(\rho )\,^t  \nabla u) + \dfrac{1}{2} \nabla( \lambda(\rho) \dive u)=0, &\label{NSAvold}
\end{eqnarray}
with 
\begin{equation}\label{NSA-constraintold}
w= u + \nu \nabla \mu(\rho)/\rho, \qquad
v=\nabla\mu(\rho)/\rho
\end{equation} 
which we call the Navier-Stokes-Korteweg augmented system in all the sequel. 

\medskip

\noindent {\it The Quantum Navier-Stokes Equations.}
Note that with the choice $K(\rho)=1/\rho$, which gives $\mu(\rho)=\rho$ and $\lambda(\rho)=0$, system \eqref{NSArhoold}-\eqref{NSA-constraintold} becomes
\begin{eqnarray}
&\partial_t\rho+\dive (\rho u)=0, & \label{NSAQrho} \\
&\partial_t(\rho \, w)+\dive (\rho \, w \otimes u)+\nabla(p(\rho))-\nu \,  \dive (\rho \nabla w)  =
 (\varepsilon^2 -\nu^2)\, \dive (\rho \nabla v), & \label{NSAQw} \\
&\partial_t(\rho  \, v)+\dive (\rho \, v \otimes u) +  \dive (\rho \,^t  \nabla u) =0, &\label{NSAQv}
\end{eqnarray}
with the constraints
\begin{equation} \label{NSAQ-constraint}
w=u+\nu \nabla\log\rho, \qquad v=\nabla\log\rho
\end{equation} which is the augmented formulation of the compressible barotropic quantum Navier-Stokes system:
\begin{eqnarray}
 &\partial_t\rho+\dive (\rho u)=0, & \label{NSQI0} \\
&\partial_t(\rho u)+\dive (\rho u \otimes u)+\nabla p(\rho)-2 \, \nu  \, \dive (\rho D(u))=2 \, \varepsilon^2 \, \rho \, \nabla\left(\dfrac{\Delta \sqrt{\rho}}{\sqrt{\rho}}\right). & \label{NSQI1}
\end{eqnarray}
In \cite{Dong10, Jiang11, Jue10}, the global existence of weak solutions to \eqref{NSQI0}-\eqref{NSQI1} has been shown following the idea introduced in \cite{BrDeLi2003} by testing the momentum equation by $\rho \, \phi$ with $\phi$ a test function. The problem of such formulation is that it requires $\gamma >3$ for $d=3$ which is not a suitable assumption for physical cases. In \cite{BDZ} the authors show the existence of solutions for \eqref{NSQI0}-\eqref{NSQI1} without quantum term ({\it i.e.} for $\varepsilon=0$) by adding a cold pressure term in the momentum equation.  The  cold pressure is  a suitable increasing function $p_c$ satisfying $\displaystyle \lim_{n \to 0} p_c(n) = + \infty$. The key element of the proof is a $\kappa$-entropy estimate. In  \cite{GV}, using the same strategy and a $\kappa$-entropy with $\kappa=1/2$, the existence of global weak solutions for \eqref{NSQI0}-\eqref{NSQI1} is proven without any extra assumption on $\gamma$ and the semi-classical limit  $\varepsilon$ tends to zero is performed. In \cite{V1}, A. Vasseur and C. Yu consider the compressible barotropic quantum Navier-Stokes equations with damping {\it i.e.} system \eqref{NSQI0}-\eqref{NSQI1} with additional terms in the right hand side of \eqref{NSQI1}: $-r_0u-r_1\rho|u|^2u$.  They  prove the global-in-time existence of weak solutions and their result is still valuable in the case $r_1=0$ . Their proof is based on a Faedo-Galerkin approximation (following the ideas of \cite{Jue10}) and a Bresch-Desjardins entropy (see \cite{BrDe2003, BrDeLi2003}). In \cite{V2}, the authors use the result obtained in  \cite{V1} and pass to the limits $\varepsilon, r_0, r_1$ tend to zero to prove the existence of global-in-time weak solutions to degenerate compressible Navier-Stokes equations. Note that to prove such a result they need uniform (with respect to $r_0, r_1$) estimates to pass to the limit $r_0, r_1$ tend to $0$. To this end they have to firstly pass to the limit $\varepsilon$ tends to $0$.  The reader interested by the compressible 
Navier-Stokes equations with density dependent viscosities is also referred to the interesting paper \cite{LiXi}.
  Recently in \cite{LaVa} and \cite{AnSp}, global existence of weak solutions for the quantum Navier-Stokes equations  \eqref{NSQI0}-\eqref{NSQI1} has been proved without drag terms and without any cold pressure. In the first paper, the method is based on the construction of weak solutions that are renormalized in the velocity variable. Note that the construction being uniform with respect to the Planck constant, the authors also perform the semi-classical limit to the associated compressible Navier-Stokes equations. Note also the recent paper \cite{AnSp} concerning the global existence for the quantum Navier-Stokes system where they use in a very nice way the mathematical structure of the equations.  It is important to remark that a global weak solutions  of the quantum Navier-Stokes equations in the sense of \cite{LaVa} is also weak solution of the augmented system (due to the regularity which is envolved allowing to write the equation on   the drift velocity $v$).   Remark also that there exists no global existence result  of weak solutions  for the compressible Navier-Stokes-Korteweg system with constant viscosities even in the two-dimensional in space case.

 \medskip
 
 \noindent {\bf Main objectives of the paper.} In this paper, to the author's point of view,  there are several interesting and new results. First starting with the global weak solutions of the quantum Navier-Stokes equations constructed in \cite{LaVa} (which is a $1/2$-entropy solution in the sense of \cite{BDZ}) we show at the viscous limit the existence of a dissipative solution for the  quantum Euler system letting the viscosity goes to zero. 
   This gives the first  global existence  result of dissipative solution for the quantum Euler system obtained from a quantum Navier-Stokes  type system. Note that in \cite{DFM}, it is proved 
the existence of infinite dissipative solutions of such inviscid quantum system. Here we
present a way to select one starting from a Navier-Stokes type system.
    Secondly, we develop relative entropy estimates for general cases of the Euler-Korteweg and the Navier-Stokes-Korteweg systems extending the augmented formulations introduced recently in \cite{BNV} and \cite{BNV2}: more general viscosities and  third order dispersive terms. 
    This gives a more simple procedure to perform relative entropy  than the one developped in  \cite{GLT, DFM}  for the Euler-Korteweg system but asks to start with 
an augmented version of the Euler-Korteweg system. This allows us to provide
a weak-strong uniqueness result for the Euler-Korteweg and Navier-Stokes-Korteweg systems.  
    
    This also helps to get rid the concavity assumption on $1/K(\rho)$ which is strongly used in \cite{GLT}. For the interested readers, we provide a comparison of the quantities appearing in our relative entropy to the ones introduced in \cite{GLT} and remark that they are equivalent under the assumptions made in  \cite{GLT}. Note that to perform our calculations for  the Navier-Stokes-Korteweg system, we need to generalize in a non-trivial way  the identity (5) in \cite{BNV}: see Proposition 30 for the generalized identity.

\medskip

For reader's convenience, let us explain the simple idea behind all the calculations. The kinetic energy corresponding to the Euler-Korteweg system reads
$$ \int_\Omega \left(\frac{1}{2} \rho |u|^2 + H(\rho) + K(\rho) |\nabla\rho|^2\right)$$
with 
$$\displaystyle H(\rho)= \rho  \int_1^\rho \dfrac{p(z)}{z^2} dz.$$
In \cite{GLT}, they consider that it is an energy written in terms of $(\rho,u,\nabla\rho)$
and they write a relative entropy playing with these unknowns. In our calculations, we write
the kinetic energy as follows
$$ \int_\Omega \left(\frac{1}{2} \rho |u|^2 + H(\rho)  + \rho |v|^2\right)$$
with $v= \sqrt{K(\rho)}\nabla\rho/\sqrt\rho$ and we consider three quantities $\rho, u$ and $v$. This motivates to write an augmented system $(\rho,u,v)$ and to modulate the energy through these three unknowns. This gives a simplest way to define an appropriate relative entropy quantity compared to \cite{GLT} and \cite{DFM} and allows to relax the concavity assumption on $1/K(\rho)$ made in the part concerning Euler-Korteweg system in \cite{GLT}.  Our result covers capillarity coefficient under the form
$$ K(\rho) \approx \rho^s \hbox{ with } s+2\le \gamma \hbox{ and } s\ge -1.$$
  Finally our result makes the link between  Euler-Korteweg system and Navier-Stokes-Korteweg system. After proving the global existence of $1/2$-entropy solutions of the general Navier-Stokes-Korteweg system (this is the subject  of a forthcoming paper \cite{BVY} still in progress: the case $K(\rho)=1/\rho$ has been recently proved in \cite{LaVa}), this could give  the  mathematical justification  of a physical dissipative solution of the Euler-Korteweg equations obtained from  $1/2$-entropy solutions of the Navier-Stokes-Korteweg equations in the spirit of \cite{BDZ}.  Note also the other
 interesting result in \cite{AnSp} on the Quantum-Navier-Stokes equations but under hypothesis between the magnitude of the viscous and capillarity coefficients.
  Let us also mention that our relative entropies could be helpful for other singular limits as explained  in the book \cite{FeNo} in the case of constant viscosities.

\bigskip
 
The paper is organized as follows.  In Section \ref{EstDef}, we provide energy estimates and the definition of weak solutions for the augmented Euler-Korteweg and Navier-Stokes-Korteweg systems. In Section \ref{sec_EuK}, we give the definition of the relative entropy formula and we established the associated estimate. This one is used to define what we call a dissipative solution for the Euler-Korteweg system and we established a weak/strong uniqueness result. The same results are obtained for the Navier-Stokes-Korteweg system in Section \ref{sec_NSK}. In Section \ref{sec_limite} we use the previous results to show the limit when the viscosity tends to zero in the quantum Navier-Stokes system.  Finally we give in Appendix some technical lemmas on modulated quantities and a comparaison between the relative entropy developed here and the one used in \cite{GLT, DFM}, and we state the definitions used for the operators.


\section{Energy estimates and definition of weak solutions.} \label{EstDef} In this subsection we give the energy equalities for the augmented Euler-Korteweg and Navier-Stokes-Korteweg systems. They will be used in the following to establish the estimates for the relative entropy associated to each one. We also define weak solutions concept  for the two augmented systems.  First of all, let us recall the definition of the function $H$ called the enthalpy by 
 $$H(\rho)= \rho e(\rho)= \rho \displaystyle  \int_1^\rho \dfrac{p(z)}{z^2} \, dz.$$ Namely we have: 
$$\rho H'(\rho)-H(\rho)=p(\rho), \, \qquad  H''(\rho)=\dfrac{p'(\rho)}{\rho}.$$ 
 To be more precise, since $p(\rho)=\rho^{\gamma}$ with $\gamma >1$, this yields to $H(\rho)=\dfrac{1}{\gamma-1}p(\rho)$.
\medskip

\noindent {\bf Euler-Korteweg system.} For the augmented Euler-Korteweg system we can show the following formal proposition.
\begin{prop} \label{NRJEK}
All strong enough solution $(\rho,u,v)$ of system \eqref{EKA1old}--\eqref{EKA4old}  satisfies:
\begin{equation*}
\dfrac{dE_{EuK}(\rho,u,v)}{dt}  =0,
\end{equation*}
where $E_{EuK}$ is the natural energy density given by
\begin{eqnarray}
\label{estentE}
E_{EuK}(t)=E_{EuK}(\rho,u,v)  
  =  \int_{\Omega} \left(\dfrac{1}{2}  \, \rho  \, |u|^2+ \dfrac{1}{2} \,  \eps^2  \, K(\rho)| \nabla \rho|^2  +H(\rho)\right).
\end{eqnarray}
\end{prop}

\proof  It suffices to take the scalar product of the equation related to $u$ by $u$ 
and the equation related to $v$ by $v$ and integrate in space using the mass equation,
the symmetry of $\nabla v$ and the relation $ \rho |v|^2 =  K(\rho) |\nabla \rho|^2.$
\cqfd

\noindent {\it Global weak solutions of the augmented system.} Assumption between $K(\rho)$ and $p(\rho)$ will be required to define global weak solutions of the augmented version of the Euler-Korteweg system namely:
$$K(\rho)= [\mu'(\rho)]^2/\rho  \hbox{ with }  \mu(\rho)= \rho^{(s+3)/2}
 \qquad \hbox{ and }\qquad  p(\rho) = \rho^\gamma$$
with
$$s+2 \le \gamma, \qquad s\ge -1 \hbox{ and } \gamma >1.$$
Assume the initial density $\rho_0$ positive and in $L^1(\Omega)$ namely  
$$ \rho_0 \ge 0  \qquad  \hbox{ and }\qquad  \int_\Omega \rho_0 < + \infty$$
and
$$E_{\rm EuK}(\rho_0,u_0,\bar{v_0})<+ \infty$$
where $\bar{v_0}$ and $u_0$ is zero where $\rho_0$ vanishes.
We can define global weak solutions of the augmented version of the Euler-Korteweg system as solutions
satisfying for a.e  $t\in [0,T]$: 
$$E_{\rm EuK}(\rho,u,\bar{v})(t) \le E_{\rm EuK}(\rho,u,\bar{v})\vert_{t=0} <+ \infty$$
with
$$\rho\ge 0  \qquad  \hbox{ and }\qquad 
    \int_\Omega \rho = \int_\Omega \rho_0
   \qquad \hbox{ and }
    \qquad  \sup_{t\in (0,T)} \int_\Omega \mu(\rho) < + \infty 
$$
 and satisfying  the following augmented system in a distribution sense
\begin{eqnarray} 
 & \partial_t\rho+\dive (\rho \, u)=0, &  \label{EKA1} \\
&   \partial_t(\rho \, u)  +\dive (\rho \, u \otimes u)+\nabla p(\rho) =
   \eps \,  {\rm div} \Bigl( \,{\mathbb T}^{EuK}(\bar{v})
  + \displaystyle \frac{\lambda(\rho)}{2\mu(\rho)}
         {\rm Tr}\, ({\mathbb T}^{EuK}(\bar{v}))\Bigr)
 & \label{EKA2} \\ 
&\partial_t(\rho  \, \bar{v})+\dive (\rho \, \bar{v} \otimes u) 
         =  - \varepsilon \, {\rm div} \Bigl( ({\mathbb T}^{EuK}(u))^t
  + \displaystyle \frac{\lambda(\rho)}{2\mu(\rho)} \, 
    {\rm Tr}\, ({\mathbb T}^{EuK}(u))\Bigr)
&  \label{EKA3} 
\end{eqnarray}
with 
\begin{equation}\label{EKA4}
 \lambda(\rho)=2(\rho \, \mu'(\rho)-\mu(\rho)), 
\qquad \bar{v} =  \varepsilon \nabla \mu(\rho)/\rho.
\end{equation}
where the tensor valued function ${\mathbb T}^{EuK}(\theta)$ (for $\theta=u$ and $\bar{v}$) is defined  through the following relation
$$ {\mathbb T}^{EuK}(\theta) 
          = \left[ \nabla(\mu(\rho)\, \theta)   -  \frac{1}{\varepsilon} \rho \theta\, \otimes \bar{v} )\right]$$
 with
 $$ {\mathbb T}^{EuK}(\theta) \in L^\infty(0,T; W^{-1,1}(\Omega)).$$

\medskip

\noindent {\it Important property.} Note that the Energy estimate provides the bound
$L^\infty(0,T;L^\gamma(\Omega))$ on $\rho$ and thus  $\mu(\rho)/\sqrt(\rho) \in L^\infty(0,T;L^2(\Omega))$ and thus using the mass quation $\mu(\rho) \in 
L^\infty(0,T; L^1(\Omega))$.

\bigskip

\noindent {\bf Navier-Stokes-Korteweg system.} Concerning the augmented Navier-Stokes-Korteweg system \eqref{NSArhoold}--\eqref{NSA-constraintold}, defining the energy 
\begin{eqnarray}
\label{entropNS}
E^{\varepsilon,\nu}_{NSK}(t) = E^{\varepsilon,\nu}_{NSK}(\rho,v,w) 
& = & \int_{\Omega} \left( \dfrac{\varepsilon^2- \nu^2}{2}  \rho \, |v|^2+\dfrac{\rho}{2}|w| ^2 +H(\rho) \right), \nonumber 
\end{eqnarray}
we have the following formal equality
\begin{prop}\label{prop_entNSK}
Let $(\rho,v,w)$  be a strong enough solution of  \eqref{NSArhoold}-\eqref{NSA-constraintold}  
we have
\begin{eqnarray*}
\dfrac{dE^{\varepsilon,\nu}_{NSK}}{dt}(\rho,v,w) + \nu \, \displaystyle \int_\Omega \left( \mu(\rho)\left(|\nabla u|^2+ (\varepsilon^2- \nu^2)|\nabla v |^2\right) +\mu'(\rho)H''(\rho) |\nabla \rho|^2 \right) \\
+\nu \, \displaystyle \int_\Omega \left(\dfrac{\lambda(\rho)}{2}\left((\dive(u))^2+(\varepsilon^2-\nu^2)(\dive(v))^2\right)\right)=0.
\end{eqnarray*}
\end{prop}
\noindent It suffices to take the scalar product of \eqref{NSAwold} with $w$ and to take the scalar product of \eqref{NSAvold} by $(\varepsilon^2-\nu^2) v$, using the expressions of $w$ and $v$,  integrate in space and sum to prove the result using the mass equation.

\bigskip

\noindent  {\it Global weak solutions of the augmented system.}
 Looking at new unknowns $(\rho,\bar{v}, w)$ with $\bar{v}= \sqrt{\varepsilon^2-\nu^2}$,
assumption between $K(\rho)$ and $p(\rho)$ will be required to define global weak solutions of the augmented version of the Navier-Korteweg system namely:
$$K(\rho)= [\mu'(\rho)]^2/\rho  \hbox{ with }  \mu(\rho)= \rho^{(s+3)/2}
 \qquad \hbox{ and }\qquad  p(\rho) = \rho^\gamma$$
with
$$s+2 \le \gamma, \qquad s\ge -1 \hbox{ and } \gamma >1.$$
Note that with this constraint on $\mu(\rho)$, we have
$$\lambda(\rho)/\mu(\rho) 
   = 2(\mu'(\rho)\rho -\mu(\rho))/\mu(\rho) 
  = (s+1) = {\rm Cst} \ge 0$$
Assume the initial density $\rho_0$ positive and in $L^1(\Omega)$ namely  
$$ \rho_0 \ge 0,  \qquad   \int_\Omega \rho_0 < + \infty $$
and
$$E_{\rm NSK}(\rho_0,\bar{v_0}, w_0)<+ \infty$$
with 
$$E_{\rm NSK} (\rho_0, \bar{v_0}, w_0)= [E_{\rm NSK} (\rho, \bar{v}, w)]_{t=0}
     = \bigl[\int_\Omega \rho |\overline{v}|^2 + \rho |w|^2 
     + H(\rho)\Bigr]_{t=0} 
       = \int_\Omega \rho_0 |\overline{v_0}|^2 + \rho_0 |w_0|^2  + H(\rho_0).$$
We can define global weak solutions of the Augmented version of the Navier-Korteweg system as solutions satisfying, for  $t\in [0,T]$,  it satisfies a.e $\tau\in [0,t]$
\begin{eqnarray} \label{energyNSK}
E_{NSK}(\rho,\bar{v},w)(\tau) 
 + \nu \, \displaystyle \int_0^t \int_\Omega \left(\left(|{\mathbb T}(w)|^2
  + |{\mathbb T}(\bar{v})|^2\right)
  + \frac{1}{\varepsilon^2-\nu^2}\frac{\rho \, p'(\rho)}{\mu'(\rho)} |\bar{v}|^2 \right)
     \nonumber \\
\hskip1cm +\nu \, \displaystyle \int_0^t\int_\Omega \left(\dfrac{\lambda(\rho)}{2\mu(\rho)}\left(|{\rm Tr}\, ({\mathbb T}(w))|^2
   +  |{\rm Tr}\, ({\mathbb T}(\bar{v}))|^2\right)\right)
\le E_{NSK}(\rho,v,w)(0)
\end{eqnarray}
where
$$E_{NSK} (\rho,\bar{v}, w) = \int_\Omega \rho |\bar{v}|^2 + \rho |w|^2 
     + H(\rho)
$$ 
$$ \rho \ge 0, \qquad  \int_\Omega \rho = \int_\Omega \rho_0 < + \infty,
     \qquad 
      \sup_{t\in (0,T)} \int_\Omega \mu(\rho) < +\infty.$$
 The augmented system in the distribution senses as follows
 \begin{eqnarray}
&\partial_t\rho+\dive (\rho u)=0, & \label{QNSNF1} \\
&\partial_t(\rho \, w)+\dive (\rho \, w \otimes u)+\nabla(p(\rho))
  -\nu \,  \dive \bigl[\sqrt {\mu(\rho )} {\mathbb T}(w)  -\dfrac{\lambda(\rho)}{2\mu(\rho)} 
      \sqrt{ \mu(\rho)}{\rm Tr} ({\mathbb T}(w)) {\rm Id}\bigr] \nonumber   \\
& =
 \sqrt{\varepsilon^2 -\nu^2}\, \dive
 \Bigl[\sqrt {\mu(\rho )} {\mathbb T}(\bar{v})  + \dfrac{\lambda(\rho)}{2\mu(\rho)} 
      \sqrt{ \mu(\rho)}{\rm Tr} ({\mathbb T}(\bar{v}))  {\rm Id}
 \Bigr], &   \label{QNSNF2} \\
&\partial_t(\rho  \, \bar{v})+\dive (\rho \, \bar{v} \otimes u) 
  - \nu  \dive \bigl[\sqrt {\mu(\rho )} T(\bar{v})  +\dfrac{\lambda(\rho)}{2\mu(\rho)} 
      \sqrt{ \mu(\rho)}{\rm Tr} \, (T(\bar{v}))  {\rm Id} \bigr]
      &  \nonumber \\
      & = - 
 \sqrt{\varepsilon^2 -\nu^2}\, 
 \dive \Bigl[\sqrt {\mu(\rho )} ({\mathbb T}(w))^t  +\dfrac{\lambda(\rho)}{2\mu(\rho)} 
      \sqrt{ \mu(\rho)}{\rm Tr} ({\mathbb T}(w))  {\rm Id}
 \Bigr], &   \label{QNSNF3} 
\end{eqnarray}
with 
\begin{equation}\label{QNSNF-constraint}
w= u + \nu \nabla \mu(\rho)/\rho, \qquad
\bar{v} = \sqrt{\varepsilon^2-\nu^2} \nabla\mu(\rho)/\rho
\end{equation} 
and where the tensor valued function $T(\theta)$ (for $\theta= w$ and $\overline v)$  satisfies $\sqrt \nu  \, T(\theta)$ is bounded in $L^2(0,T;L^2(\Omega))$  and satisfies 
the following relation
$$ \sqrt{\mu(\rho)} T(\theta) = \nabla(\mu(\rho) \, \theta) - 
     \frac{1}{\sqrt{\varepsilon^2-\nu^2}} \rho \theta\, \otimes \bar{v}
$$
and is chosen equal to zero when $\rho$ vanishes.

\medskip

\noindent {1) \it Important property.} Note that the Energy estimate provides the bound
$L^\infty(0,T;L^\gamma(\Omega))$ on $\rho$ and thus  $\mu(\rho)/\sqrt \rho \in L^\infty(0,T;L^2(\Omega))$ and thus using the mass quation $\mu(\rho) \in 
L^\infty(0,T; L^1(\Omega))$.

\medskip

\noindent {2) \it Important Remark.} Let us remark that for the global weak solutions
of the Navier-Stokes-Korteweg, the following equation is satisfied in the distribution 
sense
\begin{eqnarray} \label{mu}
\nu [\partial_t \mu(\rho) + \dive (\mu(\rho) \, u) +  
                \frac{\lambda(\rho)}{2\mu(\rho)} \sqrt{\mu(\rho)}
                {\rm Tr} \,( {\mathbb T}(u)) ] = 0
\end{eqnarray}
where $u = w - \nu \, \bar{v}/ \sqrt{\varepsilon^2-\nu^2}.$
Taking the gradient of  Equation
\eqref{mu}, we get
$$\nu [\partial_t  \nabla \mu(\rho) +  \dive({}^t\nabla(\mu(\rho) u)) 
             + \nabla \bigl(\frac{\lambda(\rho)}{2\mu(\rho)} \sqrt{\mu(\rho)} 
                     {\rm Tr}({\mathbb T}(u))
                 \bigr)]
             = 0$$
 and therefore by definition of $\sqrt{\mu(\rho)} {\mathbb T}(u)$ and expression
 of $v$, we 
 can write
 $$\nu [\partial_t( \rho v) + \dive (\rho v\otimes u)
             + \dive (\sqrt{\mu(\rho)} ({\mathbb T}(u))^t)
             + \nabla \bigl(\frac{\lambda(\rho)}{2\mu(\rho)} \sqrt{\mu(\rho)} 
                     {\rm Tr}({\mathbb T}(u))
                 \bigr)]
             = 0.$$
 This explain why a global weak solution of the Navier-Stokes-Korteweg
 system is also global weak solution of the augmented Navier-Stokes-Korteweg
 system.
 


\section{The Euler-Kortewg System : relative entropy and dissipative solution}\label{sec_EuK}

In this section, we consider the problem \eqref{EK1}--\eqref{EK2} through its augmented formulation \eqref{EKA1}--\eqref{EKA3}. The main goal of this section is to give the definition of what we call a dissipative solution for this problem. To this end we have to establish a relative entropy inequality. 


\subsection{Relative entropy inequality}

In \cite{FNL}, E. Feireisl, B.-J. Jin and A. Novotny have introduced relative entropies, suitable weak solutions and weak-strong uniqueness properties for the compressible Navier-Stokes equations with constant viscosities. The goal of this subsection is to establish a relative entropy inequality for the Euler-Korteweg System using the augmented formulation introduced in \cite{BCNV} and extending  the ideas in \cite{BNV} and \cite{BNV2} to such system
 in order to be able to  define what is called a dissipative solution.

Let us consider the following relative entropy functional, denoted $\mathcal{E}_{EuK}(\rho,u,v | r,U,V)$ and defined by
\begin{eqnarray}
&\mathcal{E}_{EuK}(t)=\mathcal{E}_{EuK}(\rho,u,v | r,U,V)(t)& \nonumber \\
&=\dfrac{1}{2} \displaystyle{\int_{\Omega} \rho  \left( |u-U|^2+\eps^2 \left|  \sqrt{\dfrac{K(\rho)}{\rho}} \nabla \rho - \sqrt{\dfrac{K(r)}{r}} \nabla r \right| ^2 \right)  \,+\int_{\Omega} H(\rho | r)}& \label{relentEK} \\
&=\dfrac{1}{2}\displaystyle{ \int_{\Omega} \rho  \left( |u-U|^2+\eps^2 \left| v - V \right| ^2 \right)  \,+\int_{\Omega} H(\rho | r),}&\nonumber
\end{eqnarray}
with 
$$ H(\rho | r) = H(\rho) - H(r) - H'(r)(\rho -r).$$
where $(\rho,u,v)$ is a weak solution of System \eqref{EKA1}--\eqref{EKA4} and
$(r,U,V)$ smooth enough target functions. Note that the definitionof the relative entropy used here is different from the one used in \cite{GLT} but we can show that the twice are equivalent in some sense for some range of the capillary coefficient. We refer to appendix \ref{annex_equivalence} for more details. Let us just say that such an energy measures the distance between a weak solution $(\rho,u,v)$ of \eqref{EKA1}--\eqref{EKA4}  to any smooth enough test function $(r,U,V)$. The goal here is to prove an inequality of type
$$\mathcal{E}_{EuK}(t)-\mathcal{E}_{EuK}(0) \leq C \int_0^t \mathcal{E}_{EuK}(\xi)d\xi, $$
with $C$ a positive constant. To this end let us first prove the following proposition.

\begin{prop}\label{relative} Let us assume that $\mu(\rho)=\rho^{(s+3)/2}$ with $\gamma \geq s+2$ and $s \geq -1$.  Let $(\rho,u,\bar{v})$ be a global weak solution to  the augmented system  \eqref{EKA1}-\eqref{EKA4}.  We have:  
\begin{eqnarray}\label{equa17}
 \displaystyle{\mathcal{E}_{EuK}(t)- \mathcal{E}_{EuK}(0)
\le  \int_0^t  \int_{\Omega}  \rho \, (U-u) \cdot  \partial_t U +  \int_0^t  \int_{\Omega} \rho (   \nabla U \, u ) \cdot (U-u) +  \int_0^t  \int_{\Omega} \rho \,  (\bar{V}- \bar{v} )  \cdot \partial_t V}& \nonumber \\
  + \displaystyle {   \, \int_0^t  \int_{\Omega} \rho (   \nabla \bar{V} \, u )\cdot (\bar{V}-\bar{v}) 
 + \varepsilon \, \int_0^t \Bigl<{\mathbb T}^{EuK}(\bar{v}) 
       + \frac{\lambda(\rho)}{2\mu(\rho)} {\rm Tr} ({\mathbb T}^{EuK}({\bar v})) {\rm Id}  ; \nabla U\Bigr> }_{W^{-1,1}(\Omega) \times W^{1,+\infty}(\Omega)} & \nonumber\\
      -  \varepsilon \, \int_0^t \Bigl<( {\mathbb T}^{EuK}(u) )^t
       + \frac{\lambda(\rho)}{2\mu(\rho)} {\rm Tr} ({\mathbb T}^{EuK}(u)) {\rm Id} ; \nabla \bar{V}  \Bigr>_{W^{-1,1}(\Omega) \times W^{1,+\infty}(\Omega)}
- \displaystyle \int_0^t  \int_{\Omega} p(\rho) \dive U \nonumber \\
-  \int_0^t  \int_{\Omega} \left[\partial_t (H'(r)) (\rho-r)+\rho \nabla H'(r) \cdot u \right], \nonumber 
 \nonumber
\end{eqnarray}
for all $t \in [0,T]$ and for all smooth test functions  $(r,U,V)$ with 
$$r \in {\mathcal C}^1([0,T] \times \bar{\Omega}), \quad  r>0,  \quad U, V \in {\mathcal C}^2([0,T]\times \bar{\Omega}).$$
\end{prop}

\proof
Thanks to the global weak solutions definition given after Proposition \ref{NRJEK} we have
\begin{eqnarray*}
\mathcal{E}_{EuK}(t)- \mathcal{E}_{EuK}(0)
 & \le &
\displaystyle \int_{\Omega} \left( \dfrac{\rho}{2} |U|^2-\rho \, u  \cdot U
 + \dfrac{1}{2} \rho \eps^2 \dfrac{K(r)}{r} |\nabla r |^2 - \rho \eps^2 \sqrt{\dfrac{K(\rho)}{\rho}} \nabla \rho \cdot \sqrt{\dfrac{K(r)}{r}} \nabla r  \right) (t) \\
 & &-\displaystyle \int_{\Omega} \left( \dfrac{\rho}{2} |U|^2-\rho \, u  \cdot U
 + \dfrac{1}{2} \rho \eps^2 \dfrac{K(r)}{r} |\nabla r |^2 - \rho \eps^2 \sqrt{\dfrac{K(\rho)}{\rho}} \nabla \rho \cdot \sqrt{\dfrac{K(r)}{r}} \nabla r  \right) (0) \\
 & &- \displaystyle \int_{\Omega} \left( H(r)+H'(r)(\rho-r)\right)(t)+\int_{\Omega}  \left( H(r)+H'(r)(\rho-r)\right)(0)
\end{eqnarray*}
{\it i.e.}
\begin{eqnarray}
\mathcal{E}_{EuK}(t)- \mathcal{E}_{EuK}(0)
& \le &
\int_0^t \int_{\Omega} \dfrac{d}{dt}\left( \dfrac{\rho}{2} |U|^2-\rho \, u  \cdot U + \dfrac{1}{2} \rho \,    \, |\bar{V}|^2 - \rho \,    \, \bar{v} \cdot \bar{V}  \right)  \nonumber \\
&&-\int_0^t \int_{\Omega} \dfrac{d}{dt}\left( H(r)+H'(r)(\rho-r)\right). \label{proofprop3_1}
\end{eqnarray}
We multiply  \eqref{EKA2} by  $U$, \eqref{EKA3} by $ \bar{V}$ and we integrate with respect to time and space.  
Writting  
$$ \partial_t(\rho \,  u \cdot U)=  \partial_t(\rho \, u)\cdot  U
+ \int_\Omega \rho \, u \cdot  \partial_t U, $$
and 
$$
 \partial_t( \rho \, v \cdot  V) =  \partial_t(\rho v) \cdot V    +   \rho \, v \cdot  \partial_t V$$
and thanks to integrations by parts, we obtain
\begin{eqnarray*}
\mathcal{E}_{EuK}(t)- \mathcal{E}_{EuK}(0)
 & \le  & \int_0^t  \int_{\Omega} \partial_t \left( \dfrac{\rho}{2} |U|^2 \right)-\int_0^t \int_{\Omega} \rho \, u \cdot  \partial_t U -\int_0^t \int_\Omega  \rho (  \nabla U \, u ) \cdot u \\ 
 & & +\eps \int_0^t
   \Bigl<  {\mathbb T}^{EuK}(v) + \frac{\lambda(\rho)}{2\mu(\rho)} 
         {\rm Tr} ({\mathbb T}^{EuK}(v)) {\rm Id}  ;  \nabla U \Bigr>_{W^{-1,1}(\Omega)\times W^{1,\infty}(\Omega)}\\
& &  + \int_0^t \int_{\Omega} \partial_t \left(  \dfrac{\rho}{2} |\bar{V}|^2 \right) - \int_0^t \int_{\Omega} \rho \, \bar{v}  \cdot\partial_t \bar{V}  \\
 & &-   \int_0^t \int_{\Omega}  \rho \left( \nabla \bar{V}  ~ u \right) \cdot \bar{v}  \\
& &   -  \eps \int_0^t
   \Bigl<  ({\mathbb T}^{EuK} (u))^t + \frac{\lambda(\rho)}{2\mu(\rho)} 
         {\rm Tr} ({\mathbb T}^{EuK}(u)) {\rm Id}  ;  \nabla \bar{V} \Bigr>_{W^{-1,1}(\Omega)\times W^{1,\infty}(\Omega)} \\
 & &  - \int_0^t \int_{\Omega} p(\rho) \dive U  -\int_0^t \int_{\Omega} \partial_t(H(r)+H'(r)(\rho-r)).
\end{eqnarray*}
Using \eqref{EKA1} and
$$\partial_t \left( \dfrac{\rho}{2} |U|^2 \right)=\dfrac{1}{2} \partial_t \rho \,  |U|^2+ \rho \,  U \cdot  \partial_t U, \,  \quad
 \partial_t \left( \dfrac{\rho}{2} |V|^2 \right)=-\dfrac{1}{2} \mbox{div} (\rho u) |V|^2+ \rho  \, V \cdot  \partial_t V,$$
thanks to integrations by parts we have
\begin{eqnarray*}
\mathcal{E}_{EuK}(t)- \mathcal{E}_{EuK}(0)
&\le & \int_0^t  \int_{\Omega}  \rho \, (U-u)  \cdot \partial_t U +\int_0^t  \int_{\Omega} \rho \,  (\bar{V}-  \bar{v} ) \cdot  \partial_t \bar{V} \\
&& +  \int_0^t  \int_{\Omega} \rho \,   ( \nabla U \,  u) \cdot  (U-u)  +  \int_0^t  \int_{\Omega} \rho  \, (\nabla \bar{V} \,  u ) \cdot (\bar{V}-\bar{v}) \\
 & & + \eps  \int_0^t 
   \Bigl< {\mathbb T}^{EuK}(\bar{v}) ; \nabla U \Bigr>_{W^{-1,1}(\Omega) \times W^{1,\infty}(\Omega)} \\
 && 
  -\eps \int_0^t  \Bigl< ({\mathbb T}^{EuK}(u))^t ; \nabla \bar{V} \Bigr>_{W^{-1,1}(\Omega) \times W^{1,\infty}(\Omega)} \\
  &&   
  +\dfrac{\eps}{2}\int_0^t  \Bigl<\frac{\lambda(\rho)}{\mu(\rho)} {\rm Tr} \,({\mathbb T}^{EuK}(\bar{v})) {\rm Id}; \nabla U
      \Bigr>_{W^{-1,1}(\Omega)\times W^{1,\infty}(\Omega)}   \\ 
 && -  \dfrac{\eps}{2}\int_0^t  \Bigl<\frac{\lambda(\rho)}{\mu(\rho)} 
   {\rm Tr} \,({\mathbb T}^{EuK}(u)) {\rm Id}; \nabla \bar{V}
      \Bigr>_{W^{-1,1}(\Omega)\times W^{1,\infty}(\Omega)}   \\ 
 &&-  \int_0^t  \int_{\Omega} p(\rho)  \,  \dive U-  \int_0^t  \int_{\Omega} \partial_t (H(r)+H'(r)(\rho-r)).  
 \end{eqnarray*}
This last inequality gives the result since with Equation \eqref{EKA1} we have:
  $$\int_{\Omega}  \partial_t (H'(r)(\rho-r))=\int_{\Omega} \left( \partial_t (H'(r)) (\rho-r)+\rho \nabla (H'(r) )\cdot u\right).$$
 \cqfd
 
 \begin{prop} \label{proposition5}
Let $(\rho,u,\bar{v})$ be a global weak solution of the augmented system \eqref{EKA1}--\eqref{EKA4} and $(r,U,\bar{V})$ be a strong solution of  
\begin{eqnarray}
&\partial_t r +\dive \, (r \, U)=0, &  \label{NSQF1} \\
&r \left( \partial_t U+   U\cdot \nabla U \right)+\nabla p(r)-\eps\left[  \dive (\mu(r) \nabla \bar{V})+ \dfrac{1}{2} \nabla( \lambda(r)  \dive \bar{V})\right]=0, &\label{NSMF2} \\
&r  \left( \partial_t \bar{V} + U \cdot  \nabla \bar{V} \right)+ \varepsilon \Bigl[ \dive (\mu(r)  \,{}^t\nabla U)+ \dfrac{1}{2} \nabla( \lambda(r)\dive U)\Bigr]=0&\label{NSMF3}
\end{eqnarray}
belonging to the class 
$$0 < {\rm inf}_{(0,T)\times \Omega}r  \le r \le {\rm sup}_{(0,T)\times \Omega} r < + \infty$$
$$    \nabla r \in L^2 (0,T;  L^\infty (\Omega)
         \cap L^1(0,T; W^{1, \infty}(\Omega))$$
$$  U   \in    L^\infty (0,T; W^{2,\infty}(\Omega))  \cap 
                      W^{1,\infty}(0,T; L^\infty(\Omega)),$$
$$ \bar{V}    \in  L^\infty(0,T; W^{2,\infty}(\Omega))
                            \cap W^{1,\infty}(0,T; L^\infty(\Omega)),$$
$$ \partial_t H'(r) \in L^1 (0,T;L^{\gamma/(\gamma-1)}(\Omega)), \qquad
       \nabla H'(r) \in L^1 (0,T; L^{2\gamma/(\gamma-1)}(\Omega))$$
and $V\vert_{t=0} = \varepsilon \nabla \mu(r_0) / r_0$.
Then we have 
\begin{eqnarray}
\mathcal{E}_{EuK}(t)- \mathcal{E}_{EuK}(0)
&\le&    \int_0^t  \int_{\Omega}  \rho(u-U) \cdot (\nabla U (U-u)) -   \int_0^t  \int_{\Omega}  \rho(\bar{V}-\bar{v}) \cdot(\nabla \bar{V} (U-u)) \nonumber \\
&{}& - \int_0^t  \int_{\Omega} (p(\rho)-p(r)-(\rho-r)p'(r))\dive U \nonumber \\
&{}&- \int_0^t  \int_{\Omega} \rho \,   (  \bar{v}-\bar{V} )  \cdot  \nabla  U   (\bar{v}-\bar{V})   +    \int_0^t  \int_{\Omega} \rho  (\bar{v}-\bar{V} )  \cdot   \nabla \bar{V} (u-U) \nonumber \\
&{}&-   \varepsilon\int_0^t  \int_{\Omega}  \rho \left(   \mu''(\rho) \nabla \rho  -  \mu''(r)\nabla(r) \right) \cdot \left(  (\bar{v}-\bar{V}) \, \dive U + (U-u) \, \dive \bar{V} \right)  \nonumber \\
&{}&- \varepsilon\int_0^t  \int_{\Omega} \rho  \,  (  \mu'(\rho)  -  \mu'(r)) \left(  (\bar{v}-\bar{V}) \, \cdot \nabla(\dive U)  +(U-u) \, \cdot \nabla(\dive \bar{V}) \right). \nonumber 
\end{eqnarray}
\end{prop}

\proof
First remark that due to the initial condition hypothesis and the regularity hypothesis on $U$, we can prove that $\bar{V}= \varepsilon \nabla \mu(r)/r.$
Multiplying  \eqref{NSMF2} by $\dfrac{\rho}{r} \, (U-u)$ and  \eqref{NSMF3} by  $ \, \dfrac{\rho}{r} \, ({\bar V} - \bar{v})$ and integrating with respect to time and space we have:
\begin{eqnarray}
\mathcal{E}_{EuK}(t)- \mathcal{E}_{EuK}(0)
 & \le &     \displaystyle -  \int_0^t  \int_{\Omega}  \rho \, (\nabla U \, (U-u)) \cdot  (U-u) 
 \nonumber \\
& &   \displaystyle -   \int_0^t  \int_{\Omega}  \rho \, (\nabla \bar{V} \, (U-u)) \cdot  (\bar{V}-\bar{v}) \displaystyle \nonumber \\
 & & + \eps  \,( I_1^{EuK}+I_2^{EuK}) +I_3^{EuK}  \nonumber 
\end{eqnarray}
with 
\begin{eqnarray*}
I_1^{EuK} & = &   \int_0^t  \int_{\Omega}  \left( \dfrac{\rho}{r} \dive (\mu(r) \nabla \bar{V})\cdot  (U-u)-   \dfrac{\rho}{r} \dive (\mu(r) {}^t\nabla U)\cdot (\bar{V}-\bar{v}) \right)  \\
& & +  \int_0^t   
\Bigl< {\mathbb T}^{EuK}(\bar{v});
           \nabla U\Bigr>_{W^{-1,1}(\Omega)\times W^{1,\infty}(\Omega)}  \\      
& &    - \int_0^t 
 \Bigl<({\mathbb T}^{EuK}(u))^t;
           \nabla \bar{V}\Bigr>_{W^{-1,1}(\Omega)\times W^{1,\infty}(\Omega)} 
\end{eqnarray*}
 \begin{eqnarray*}
2 \, I_2^{EuK}&= & \int_0^t  \int_{\Omega}  \dfrac{\rho}{r} (U-u) \cdot \nabla \left(\lambda(r) \,  \dive \bar{V} \right)  -\int_0^t  \int_{\Omega}  \dfrac{\rho}{r} (\bar{V}-\bar{v}) \cdot \nabla \left(\lambda(r)  \,   \dive U \right) \\
&{}& +  \int_0^t  \Bigl<\frac{\lambda(\rho)}{\mu(\rho)} {\rm Tr}({\mathbb T}^{EuK}(\bar{v})) {\rm Id} ;
           \nabla U\Bigr>_{W^{-1,1}(\Omega)\times W^{1,\infty}(\Omega)} \\
&&  -   \int_0^t  \Bigl<\frac{\lambda(\rho)}{\mu(\rho)}{\rm Tr}({\mathbb T}^{EuK}(u)) {\rm Id} ;
           \nabla \bar{V}\Bigr>_{W^{-1,1}(\Omega)\times W^{1,\infty}(\Omega)}
\end{eqnarray*}
\begin{eqnarray*}
I_3^{EuK}  & = &    \displaystyle  \int_0^t  \int_{\Omega}\left( - p(\rho)   \, \dive U  -   \dfrac{\rho}{r} \,  \nabla p(r) \cdot (U-u)
-    \partial_t (H'(r)) (\rho-r)- \rho \,  \nabla (H'(r)) \cdot  u \right).
\end{eqnarray*}
Using $r H''(r)=p'(r)$, we have 
 $$\dfrac{\rho}{r} \nabla p(r) =\rho \nabla (H'(r)).$$
Mutiplying  \eqref{NSQF1} by  $H''(r)$ and using $r H''(r)=p'(r)$ we obtain
 $$\partial_t H'(r)   +  \nabla H'(r) \cdot U  +  p'(r)\dive U=0.$$
Using $r H''(r)=p'(r)$  and an integration by parts, we have
 $$\int_0^t  \int_{\Omega} r \nabla H'(r) U =-  \int_0^t  \int_{\Omega} p(r) \dive U.$$
 Then, 
$$I_3^{EuK}= \int_0^t  \int_{\Omega} (p(r)-p(\rho)-(r-\rho)p'(r))\dive U.$$
We have
\begin{eqnarray} \label{I1}
I_1^{EuK} =  I_4^{EuK}+I_5^{EuK},
\end{eqnarray}
where
\begin{eqnarray*}
\varepsilon I_4^{EuK} & = &  
  \varepsilon \int_0^t  \int_\Omega \frac{\rho}{r} \mu(r)
   \Bigl[ \Delta \bar{V}  \cdot (U-u)  -  \nabla {\rm div} U \cdot (\bar{V}-\bar{v})\Bigr]\\
&&   \int_0^t  \int_{\Omega} 
    \rho (\bar{V}\cdot \nabla \bar{V}) \cdot (U-u)
    - \rho ((\bar{V}- \bar{v}) \cdot \nabla U) \cdot \bar{V} ,   \\
  \end{eqnarray*}
and using the symmetry of $\nabla \bar{v}$ and $\nabla \bar{V}$ and the definition the tensor value function ${\mathbb T}^{EuK} (u)$ and ${\mathbb T}^{EuK}(\bar{v})$ which
may be also written for $U$ and $\bar{V}$ (recalling that $\mu(\rho) \in L^\infty(0,T;L^1(\Omega))$),
\begin{eqnarray*}
\varepsilon I_5^{EuK}  & =     &   \varepsilon \int_0^t   
\Bigl< {\mathbb T}^{EuK}(\bar{v});
           \nabla U\Bigr>_{W^{-1,1}(\Omega)\times W^{1,\infty}(\Omega)}  
             - \varepsilon  \int_0^t 
 \Bigl<({\mathbb T}^{EuK}(u))^t;
           \nabla \bar{V}\Bigr>_{W^{-1,1}(\Omega)\times W^{1,\infty}(\Omega)} \\
 & = &   \varepsilon \int_0^t   
\Bigl<({\mathbb T}^{EuK}(\bar{v}))^t ;
           \nabla U\Bigr>_{W^{-1,1}(\Omega)\times W^{1,\infty}(\Omega)}  
             - \varepsilon  \int_0^t 
 \Bigl< {\mathbb T}^{EuK}(u);
           \nabla \bar{V}\Bigr>_{W^{-1,1}(\Omega)\times W^{1,\infty}(\Omega)} \\
& &  -   \varepsilon \int_0^t   
\Bigl< (({\mathbb T}^{EuK}(\bar{V}))^t;
           \nabla U\Bigr>_{W^{-1,1}(\Omega)\times W^{1,\infty}(\Omega)}  
         +  \varepsilon  \int_0^t 
 \Bigl< {\mathbb T}^{EuK}(U);
           \nabla \bar{V}\Bigr>_{W^{-1,1}(\Omega)\times W^{1,\infty}(\Omega)} \\
 &= & \varepsilon \int_0^t \int_\Omega \mu(\rho) \Bigl[(\bar{V}- \bar{v})\cdot \nabla {\rm div} U
              + (u-U) \Delta \bar{V} \Bigr]     \\
 && + \int_0^t \rho((\bar{V}- \bar{v}) \cdot \nabla U)\cdot \bar{v} -
                      \rho (\bar{v} \cdot \nabla \bar{V})\cdot (U-u)
\end{eqnarray*}
  Then we get
  \begin{eqnarray*}
\varepsilon I_1^{EuK} & = &  
   \varepsilon   \int_0^t  \int_{\Omega}  \rho  \left(\dfrac{\mu(\rho)}{\rho} -   \dfrac{\mu(r)}{r} \right) (
   (\bar{v}-\bar{V})\cdot \dive(  {}^t \nabla U )    +  (u-U) \cdot \dive (\nabla  \bar{V}) )   \\
  & &  -  \int_0^t  \int_{\Omega} \rho  (\bar{v}-\bar{V} )  \cdot   \nabla  U   (\bar{v}-\bar{V})   + \int_0^t  \int_{\Omega} \rho  (\bar{v} - \bar{V} )  \cdot   \nabla \bar{V} (u-U). 
  \end{eqnarray*}
Let us now look at $I_2^{EuK}$, we have
 \begin{eqnarray*}
2  I_2^{EuK} & =  &  \int_0^t  \int_{\Omega}    \dfrac{\rho}{r} (U-u) \cdot\nabla (\lambda (r) \dive \bar{V})-   \int_0^t  \int_{\Omega}  \dfrac{\rho}{r} (\bar{V}-\bar{v}) \cdot\nabla (\lambda (r)   \dive U) \\
 &  &  + \int_0^t  \Bigl<\frac{\lambda(\rho)}{\mu(\rho)}
       {\rm Tr}({\mathbb T}^{EuK}(\bar{v})) {\rm Id} 
         ; \nabla U \Bigr>_{W^{-1,1}(\Omega)\times W^{1,\infty}(\Omega)} \\
 & & - + \int_0^t  \Bigl<\frac{\lambda(\rho)}{\mu(\rho)}
        {\rm Tr}({\mathbb T}^{EuK}(u)) {\rm Id} 
         ; \nabla \bar{V}  \Bigr>_{W^{-1,1}(\Omega)\times W^{1,\infty}(\Omega)} \\
\end{eqnarray*}
and therefore recalling that $\lambda'(\rho) = 2 \rho \mu''(\rho)$ and playing
as for $I_1^{EuK}$, we get 
 \begin{eqnarray*}
2  I_2^{EuK} & =  & -2 \,  \int_0^t  \int_{\Omega}  \rho \left(   \mu''(\rho) \nabla \rho  -  \mu''(r)\nabla r \right) \cdot  \left(  (\bar{v}- \bar{V}) \, \dive U + (U-u) \, \dive \bar{V} \right)  \\
&  &- \int_0^t  \int_{\Omega}  \left( \lambda (\rho) -\dfrac{\rho}{r}\lambda (r)\right) 
\left(  (\bar{v}- \bar{V}) \cdot \nabla(  \dive U)  +(U-u) \cdot \nabla ( \dive \bar{V} )\right)  \\
  \end{eqnarray*}
  and therefore because $\lambda(\rho) = 2 (\mu'(\rho)\rho-\mu(\rho))$, we get
  \begin{eqnarray*}
 \varepsilon( I_1^{EuK}+I_2^{EuK}) 
 & = & - \int_0^t  \int_{\Omega} \rho ( \bar{v}- \bar{V} )  \cdot  \nabla  U (\bar{v}-\bar{V})   +  \int_0^t  \int_{\Omega} \rho  (\bar{v}-\bar{V})   \cdot   \nabla \bar{V} (u-U)  \\
   & & - \varepsilon \int_0^t  \int_{\Omega}  \rho \left(   \mu''(\rho) \nabla \rho  -  \mu''(r)\nabla(r) \right) \cdot \left(  (\bar{v}-\bar{V}) \, \dive U + (U-u) \, \dive \bar{V} \right)  \\
&  & - \varepsilon \int_0^t  \int_{\Omega} \rho  (  \mu'(\rho)  -  \mu'(r)) \left(  (\bar{v}-\bar{V}) \cdot \nabla ( \dive U ) +(U-u) \cdot \nabla ( \dive \bar{V} )\right). \\
 \end{eqnarray*} 
 This concludes the proof.
\cqfd

\begin{theorem} \label{th_EuK_1}
Let us assume $\mu(\rho)=\rho^{(s+3)/2}$ with $ \gamma \geq s+2$ and $s \geq -1$. Let $(\rho,u,v)$ be a global weak solution of the augmented system \eqref{EKA1}--\eqref{EKA4}  and $(r,U,V)$ be a strong solution of \eqref{NSQF1}--\eqref{NSMF3} in the sense of
Proposition \ref{proposition5}.
We have 
\begin{eqnarray*}
\mathcal{E}_{EuK}(t)- \mathcal{E}_{EuK}(0) \leq C(r,U,V) \int_0^t  \int_{\Omega} \mathcal{E}_{EuK}(\xi)d\xi ,
\end{eqnarray*}
where $C(r,U,V)$ is a uniformly bounded constant on $\R^+ \times \Omega$.
\end{theorem}

\noindent Using the Gronwall's Lemma, we directly obtain:
\begin{cor} \label{corenergie}
Let us assume $\mu(\rho)=\rho^{(s+3)/2}$ with $ \gamma \geq s+2$ and $s \geq -1$. 
 Let  $(\rho,u,v)$ be a global weak solution of  \eqref{EKA1}-\eqref{EKA4} and  $(r,U,V)$ a strong solution of  \eqref{NSQF1}-\eqref{NSMF3} in the sense of proposition \ref{proposition5}. Then 
$$\mathcal{E}_{EuK}(t)\leq \mathcal{E}_{EuK}(0)\exp(C \,  t),$$
with  $C=C(r,U,V)$ a uniformly bounded constant on $\R^+ \times \Omega$.
It the initial conditions coincide for $(\rho,u,v)$ and $(r,U,V)$ then $\rho=r$, 
$u=U$, $v=V$.
\end{cor}

Note that theorem \ref{th_EuK_1} is a direct consequence of proposition \ref{proposition5} and the following lemma.

\begin{lem} \label{lem6}
We assume that $\mu(\rho)=\rho^{(s+3)/2}$ with $s \geq -1$. 
Let $(\rho,u,v)$ be a global weak solution of \eqref{EKA1}--\eqref{EKA4} 
and $(r,U,V)$ be a strong solution of \eqref{NSQF1}--\eqref{NSMF3} in the
sense of Proposition \ref{proposition5} .
Then
\begin{eqnarray*}
&\displaystyle
\varepsilon \left| \int_0^t  \int_{\Omega}  \rho \left(   \mu''(\rho) \nabla \rho  -  \mu''(r)\nabla(r) \right) \cdot \left(  (\bar{v}-\bar{V}) \, \dive U + (U-u) \, \dive \bar{V} \right) \right|& \\
&\displaystyle\leq C \, \dfrac{s+1}{2} \int_0^t  \int_{\Omega} \rho (|\bar{v}-\bar{V}|^2 +|u-U|^2),&
\end{eqnarray*}
and, if $\gamma \geq 2+s$, we have
$$ \left |\int_0^t  \int_{\Omega} \rho  (  \mu'(\rho)  -  \mu'(r)) \left(  (\bar{v}-\bar{V}) \, \cdot \nabla ( \dive U ) +(U-u) \, \cdot \nabla (\dive \bar{V} )\right) \right|$$
$$\leq  C  \int_0^t   \int_\Omega  \left(H(\rho | r) +   \rho (|\bar{v}-\bar{V}|^2 +|u-U|^2) \right),
$$
where $C=C(r,U,V)$ is a uniformly bounded constant on $\R^+ \times \Omega$.
\end{lem}

\proof As $(r,U,\bar{V})$ is a strong solution of \eqref{NSQF1}--\eqref{NSMF3} 
then we can prove that $\bar{V}=\varepsilon \nabla(\mu(r))/r$. Since $\mu(\rho)=\rho^{(s+3)/2}$ and $\bar{v}=\varepsilon \nabla(\mu(\rho))/\rho$ we have   
\begin{eqnarray*}
 \varepsilon(\mu''(\rho) \nabla \rho  -  \mu''(r)\nabla r)  & = &  \frac{s+1}{2} (\bar{v}-\bar{V}),
\end{eqnarray*}
 which  gives the first part of the lemma using Young inequality.
For the second one, using Young inequality, we have:
\begin{eqnarray*}
&{}&\int_0^t  \int_{\Omega} \left| \rho  (  \mu'(\rho)  -  \mu'(r)) \left(  (\bar{v}-\bar{V}) \, \cdot \nabla ( \dive U ) +(U-u) \, \cdot \nabla (\dive \bar{V}) \right) \right| \\
&{}& \leq  C\left(\dfrac{1}{2} \int_0^t \int_\Omega \rho | \mu'(\rho)-\mu'(r) |^2 +  \int_0^t\int_\Omega \rho |\bar{v}-\bar{V} |^2 +  \int_0^t\int_\Omega \rho | u-U |^2 \right),
\end{eqnarray*}
with $C=C(U,V)$ a uniformly bounded constant on $\R^+ \times \Omega$. Using Lemma \ref{lem_tech4} in the first integral, we obtain the result.
\cqfd

Let us now give a weak-strong uniqueness result based on solutions that has been
already constructed in \cite{AM12}-\cite{CaDaSa} and \cite{BDD2007}.
\begin{theorem} \label{th_QEK} Let $(r_0, u_0) \in H^{s+1}(\Omega) \times H^s(\Omega)$
with $s>2+d/2$ with $r_0>0$ such that ${\rm curl}(r_0 u_0) = 0$.
Let $(\rho,u)$ be a global weak solution in $(0,T)\times \Omega$ of the 
Quantum-Euler system 
\begin{eqnarray}
& &  \partial_t \rho + {\rm div} (\rho u) = 0 \\
& & \partial_t (\rho u ) + {\rm div}(\rho u\otimes u) 
         = \varepsilon^2 \dive\,  (\nabla \nabla \rho 
              -  \rho \nabla\log \rho \otimes \nabla \log \rho)  
\end{eqnarray}b
 corresponding to the initial data $(r_0,r_0u_0)$ and let $(r,U)$ a local strong solution
 in $(0,T^*) \times \Omega$ of this system for the same initial data with
 $$0 < c \le r \le c^{-1} < +\infty$$
 where $c$ is a constant and 
 $$r \in {\mathcal C} ([0,T]\times H^{s+1}(\Omega)) \times 
              {\mathcal C}^1 ([0,T]\times H^{s-1}(\Omega))$$
 $$ U \in {\mathcal C} ([0,T]\times H^{s}(\Omega)) \times 
               {\mathcal C}^1 ([0,T]\times H^{s-2} (\Omega))$$
 then $\rho=r$, $u=U$ and $w=W$ on $(0,\min{\{T,T^*\}})\times \Omega$.
\end{theorem}

\noindent {\it Proof.} Let us first remark that such existence of local strong solution
has been proved for instance in \cite{BDD2007} in the whole space without the constraint
on ${\rm curl} (\rho_0 u_0)= 0$ but may be considered in the periodic case. The global existence of weak solution for the Quantum--Euler System with the constraint ${\rm curl} (\rho_0 u_0)=0$ has been obtained in two papers namely \cite{AM12} and \cite{CaDaSa}.  For a strong solution, it is not difficult to prove that it also satisfies the augmented system. Concerning the global weak solution,   it  suffices to recall the important remark given in the introduction.  Differentiating in space the mass equation in ${\mathcal D}'((0,T)\times \Omega)$  we get
$$\partial_t \nabla \rho + \nabla \dive (\rho u) = 
     \partial_t \nabla \rho + \dive ({}^t\nabla (\rho u))
       = \varepsilon^2  \dive [ \Delta \rho -  \nabla \sqrt \rho \otimes \sqrt \rho] $$
which may be written
$$\partial_t \nabla \rho + \dive (\rho \nabla \log\rho \otimes u)
                                         + \dive\bigl({}^t\nabla(\rho u) - \rho \nabla \log\rho \otimes u)=0$$
 and therefore  
 $$ \partial_t \nabla \rho + \dive (\rho \nabla \log\rho \otimes u) 
       + \dive\bigl({\mathbb T}(u)^t)=0.$$       
Using the definition $\rho \bar{v} =\varepsilon  \nabla \rho$, we can rewrite the
Quantum-Euler system and the previous relation  in its augmented form
\begin{eqnarray}
& &  \partial_t \rho + {\rm div} (\rho u) = 0 \\
& & \partial_t (\rho u ) + {\rm div}(\rho u\otimes u) 
         = \varepsilon \dive\,  {\mathbb T}(\bar{v})  \\
& &  \partial_t \rho \bar{v}  + \dive (\rho \bar{v} \otimes u) 
       + \varepsilon  \dive\bigl({\mathbb T}(u)^t)=0
\end{eqnarray}
which is the augmented version of the Quantum-Euler equations.
Thus a global weak solution of the Quantum-Euler system is a global weak
solution of the augmented Quantum-Euler system and therefore the weak-strong uniqueness corollary \ref{corenergie} may be applied due to the regularity of the
strong solution.


\subsection{Dissipative solutions and weak-strong uniqueness result}

In this subsection, we  give the definition of what we call a dissipative solution for the Euler-Korteweg System. We recall that $\mathcal{E}_{EuK}(t)$ stands for
$$\mathcal{E}_{EuK}(t) =\mathcal{E}_{EuK}(\rho,u,v |  r,U,V)(t)$$
defined in \eqref{relentEK}.
Let $U$ be a smooth function, then we solve the transport equation for $r$ for  the initial
data $r_0$ such that $0<r_0<+\infty$. We then define
 the function $\mathscr{E}$ as
\begin{equation}
 \mathscr{E}(r,U)= 
 r \left( \partial_t U+ U \cdot  \nabla U \right)+\nabla p(r)-\eps^2  \dive (\mu(r) \, ^t  \nabla V)+\dfrac{1}{2} \nabla(\lambda(r) \dive V),  \label{soldisI2EuK} 
\end{equation}
with $r\, V=\nabla(\mu(r))$. Then we can prove  differentiating \eqref{NSQF1}, that
\begin{equation}
0= r  \left( \partial_t V+ U \cdot \nabla V  \right)+ \dive (\mu(r)  \,{}^t\nabla U)+ \dfrac{1}{2} \nabla(\lambda(r)  \dive U). \label{soldisI3EuK} 
\end{equation}

\begin{defi}\label{Defi10EK}
   Let us assume $\mu(\rho)=\rho^{(s+3)/2}$ ({\it i.e.} $K(\rho)= \frac{(s+3)^2}{4} \rho^s$)  with $ \gamma \geq s+2$ and $s \geq -1$. 
 Let $\rho_0$ and $u_0$ smooth enough. The pair $(\rho,u,v)$   is a dissipative solution of the Euler-Korteweg System corresponding to the initial conditions 
$$\rho \vert_{t=0} = \rho_0, \qquad 
    \rho u \vert_{t=0} = \rho_0 u_0, \qquad
    \rho v\vert_{t=0} = \sqrt{\rho_0 K(\rho_0)}\nabla \rho_0.
$$
if the triplet $(\rho,u,v)$ satisfies 
$$\mathcal{E}_{EuK}(t) \leq    \mathcal{E}_{EuK}(0) \exp(C \, t) +b_{EuK}(t)+C \,  \int_0^t b_{EuK}(\xi) \exp( C \, (t-\xi) )\,d\xi, $$
with $C=C(\eps^2,r,U,V)$ a uniformly bounded constant on $\R^+ \times \Omega$,
and where
$$b_{EuK}(t)= \int_0^t  \int_{\Omega} \dfrac{\rho}{r}|\,  \mathscr{E} \,\cdot\,  (U-u) |,$$ 
for all strong enough $U$ test function  and $(r,\mathscr{E})$ given respectively through  \eqref{NSQF1} and  \eqref{soldisI2EuK} and  the identity \eqref{soldisI3EuK}.
\end{defi}

 As a direct consequence, we can establish the following weak-strong uniqueness property (see \cite{G}).
\begin{theorem}
Let us assume $\mu(\rho)=\rho^{(s+3)/2}$ ({\it i.e.} $K(\rho)=\dfrac{(s+3)^2}{4}\rho^s$)  with $ \gamma \geq s+2$ and $s \geq -1$. 
Let us consider a dissipative solution  $(\rho,u,v)$  to the Euler-Korteweg system \
satisfying the initial conditions 
$$\rho \vert_{t=0} = \rho_0, \qquad 
    \rho u \vert_{t=0} = \rho_0 u_0, \qquad
    \rho v\vert_{t=0} = \sqrt{\rho_0 K(\rho_0)}\nabla \rho_0.
$$
Let us assume that $(r,U)$ is a strong solution of  \eqref{NSQF1} and 
\begin{equation}
\label{EKWSU1}
r \left( \partial_t U+ \nabla U \,U \right)+\nabla p(r)- \eps^2  r \nabla\left(K(r)\Delta r+\dfrac{1}{2}K'(r)|\nabla r|^2\right)=0
\end{equation}
with the regularity given in proposition \ref{proposition5} where we denote 
$\bar{V}=\varepsilon \nabla(\mu(r))/r$ and with $(\rho_0,u_0) \in W^{2,\infty}(\Omega)\times W^{1,\infty}(\Omega)$.
   If $r\vert_{t=0} =\rho_0$, $U\vert_{t=0} =u_0$ then  $\rho=r, \, u=U$ and $v=V$, which means that the problem satisfies  a dissipative-strong uniqueness property.
\end{theorem}

\proof 
If $(r,U)$ is a strong solution of \eqref{NSQF1}, \eqref{EKWSU1} then $\mathscr{E}=0$ and $b_{EuK}(t)=0$. We have
\begin{equation}
\label{proffwsuniqEQ}
0 \leq \mathcal{E}_{EuK}(t) \leq    \mathcal{E}_{EuK}(0) \exp(C \, t).
\end{equation}
If $r(t=0)=\rho_0$, $U(t=0)=u_0$ then $v(t=0)=V(t=0)$ and $\mathcal{E}_{EuK}(0)=0$, then this leads to $\rho=r$, $u=U, \, v=V$ using \eqref{proffwsuniqEQ}.
\cqfd

Note that,  as already mentioned before, all the results and definitions of this section are still valid for the compressible quantum Euler System. Indeed it corresponds to the special case $K(\rho)=1/\rho$ in the Euler-Korteweg System for which the assumption $2+s \leq \gamma$ is satisfied since $s=-1$ and $\gamma >1$. In particular we have the following definition of what we call a dissipative solution of the quantum Euler system. This one will be used in section \ref{sec_limite}.
\begin{defi}\label{Defi10}
Let $\rho_0$ and $u_0$ smooth enough. The triple $(\rho,u,v)$   is a dissipative solution of the quantum Euler system \eqref{EQ1}-\eqref{EQ2} corresponding to the initial conditions $$\rho \vert_{t=0} = \rho_0, \qquad 
    \rho u \vert_{t=0} = \rho_0 u_0, \qquad
    \rho v\vert_{t=0} = \rho_0 \nabla \log \rho_0.
$$
if the triplet $(\rho,u,v)$ satisfies 
$$\mathcal{E}_{EuQ}(t) \leq    \mathcal{E}_{EuQ}(0) \exp(C \, t) +b_{EuQ}(t)+C \,  \int_0^t b_{EuQ}(\xi) \exp( C \, (t-\xi) )\,d\xi, $$
where  $\mathcal{E}_{EuQ}(0)=  \mathcal{E}_{EuQ}\vert_{t=0}$ and
with a constant $C=C(\eps^2,r,U,V)$ uniformly bounded on $\R^+\times \Omega$, and 
$$\mathcal{E}_{EuQ}(t)=\mathcal{E}_{EuK}(t), \hbox{ for } K(\rho)=1/\rho,$$
$$b_{EuQ}(t)= \int_0^t  \int_{\Omega} \dfrac{\rho}{r}|\,  \mathscr{E} \,\cdot\,  (U-u) |,$$ 
for all smooth $U$ and $(r,V,\mathscr{E})$ defined respectively through  \eqref{NSQF1} and 
\begin{eqnarray}
& V= \nabla \log r, \\
&  \mathscr{E}(r,U)= \left( \partial_t U+ U \cdot  \nabla U  \right)+\nabla p(r)- \eps^2  \dive \, (r \nabla V),  \label{soldisI2} 
\end{eqnarray}
\end{defi}

\noindent {\it Remark.}
   Note that, in the definition above, since $U$ is regular and also $r$, we have $V$ which satisfies
\begin{equation} \label{soldisI3}
r  \left( \partial_t V+ U  \cdot \nabla V \right)+  \dive \,  (r  \, ^t  \nabla U)= 0.
\end{equation}



\section{The Navier-Stokes-Korteweg System: relative entropy and dissipative solution}\label{sec_NSK}

The goal of this section is  to define what we call a dissipative solution for the Navier-Stokes-Korteweg  System. To this end, we consider the augmented System \eqref{NSArhoold}-\eqref{NSAvold} and we establish a relative entropy estimate. Here the viscous term adds some difficulties compare to the case of the Euler-Korteweg system.


\subsection{Relative entropy inequality}

In this section, we establish a relative entropy inequality for a weak solution $(\rho,\bar{v},w)$ of  the augmented System \eqref{QNSNF1}-\eqref{QNSNF3}.
This will then be used to give the definition of what is called a dissipative solution for the Navier-Stokes-Korteweg system.
We define the following relative entropy functional
\begin{eqnarray*}
\mathcal{E}_{NSK}(t) & = & \mathcal{E}_{NSK}(\rho,\bar{v},w |  r, \bar{V},W) \\
&= & \dfrac{1}{2} \int_{\Omega} 
\rho \, \left(  |\bar{v}-\bar{V}|^2 +  |w-W|^2  \right)  
+ \int_{\Omega} H(\rho | r) \\
 & & +\nu \, \int_0^t \,  \int_{\Omega} \mu( \rho) \, \left(|\frac{{\mathbb T}(\bar{v})}{\sqrt{\mu(\rho)}} - \nabla \bar{V}|^2 +  |\frac{{\mathbb T}(\bar{w})}{\sqrt{\mu(\rho)}}
       -\nabla W |^2\right) \\
  & &+ \dfrac{\nu}{2}   \int_0^t \,  \int_{\Omega}  \lambda(\rho) 
       \left( (\frac{{\rm Tr}\,{\mathbb T}(\bar{v})}{\sqrt{\mu(\rho)}} -\dive \bar{V} )^2 +(\frac{{\rm Tr}\,{\mathbb T}(w)}{\sqrt{\mu(\rho)}} -\dive W)^2 \right). 
  \end{eqnarray*}

\begin{prop} \label{Prop19}
Any global weak solution $(\rho,\bar{v},w)$ of the augmented system \eqref{QNSNF1}-\eqref{QNSNF-constraint} satisfies the follo\-wing ine\-quality
for all $t \in [0,T]$ and for any test functions $$r \in C^1([0,T] \times \bar{\Omega}), \quad  r>0, \quad \bar{V}, \, W \in C^2([0,T] \times \bar{\Omega}),$$
 \begin{eqnarray} 
 \mathcal{E}_{NSK}(\rho,\bar{v},w |  r, \bar{V},W)(t)  
 &  \leq & 
  \mathcal{E}_{NSK}(\rho,\bar{v},w |  r, \bar{V},W)(0) \nonumber \\
  & {}  & +\displaystyle  \int_0^t  \int_{\Omega}  \rho \left(  \partial_t \bar{V} \cdot (\bar{V}-\bar{v}) 
 +   (\nabla \bar{V} \, u) \cdot  (\bar{V}-\bar{v})   \right)  \nonumber  \\
  & {}  &   + \displaystyle    \int_0^t  \int_{\Omega} \rho \left(  \partial_t W \cdot (W- w ) + (\nabla W \, u) \cdot (W-w) \right)  \nonumber  \\
 & {}  &  + \nu \displaystyle   \int_0^t  \int_{\Omega}  \mu(\rho)  \,  \left(  |\nabla \bar{V}|^2 +|\nabla W|^2\right) -  \sqrt{\mu(\rho)}  \left( {\mathbb T}(\bar{v}) : \nabla \bar{V} 
 + {\mathbb T} (w)  : \nabla W \right)   \nonumber   
  \end{eqnarray}
 \begin{eqnarray}
 & {}  & +  \sqrt{\eps^2 - \nu^2} \displaystyle   \int_0^t \int_{\Omega}
   \sqrt{ \mu(\rho)} \left({\mathbb T}(\bar{v}) :  \nabla W 
        -  ({\mathbb T}(\bar{w}))^t  : \nabla \bar{V}
 \right)   \nonumber  \\
  & {}  & +  \sqrt{\eps^2 - \nu^2}  \displaystyle   \int_0^t \int_{\Omega}
\dfrac{\lambda(\rho)}{2\sqrt{\mu(\rho)}}  \left({\rm Tr}({\mathbb T}(\bar{v})\dive W
    - {\rm Tr}({\mathbb T}(w)) \dive \bar{V}  \right)  \nonumber  \\
   & {}   & +\dfrac{\nu}{2}  \displaystyle   \int_0^t \int_{\Omega}    \lambda(\rho) \left( (\dive \bar{V})^2 + ( \dive W )^2\right)  \\
   & {} & - \dfrac{\nu}{2}  \displaystyle   \int_0^t \int_{\Omega}  \frac{\lambda(\rho)}{\sqrt{\mu(\rho)}} \left({\rm Tr}({\mathbb T}(\bar{v}) \dive\bar{V} + {\rm Tr}({\mathbb T}(w))\dive W \right)  \nonumber  \\
  & {}   &  -  \displaystyle   \int_0^t  \int_{\Omega}  \left( \partial_t( H'(r)) (\rho-r)+ \rho \, \nabla (H'(r)) \cdot  \, u +p(\rho) \dive W \right)  \nonumber  \\
 & {}   &-  \nu \displaystyle   \int_0^t  \int_{\Omega}\mu'(\rho) H''(\rho) |\nabla \rho |^2.  \nonumber   
  \end{eqnarray}
\end{prop}

\noindent {\it Remark.} Note that each quantities are defined in the usual sense for weak solution $(\rho,\bar{v},w)$  and
regular test functions $(r,\bar{V},W)$ as chosen in the proposition above. The main difference compared to the Euler-Korteweg
system is that here we control $\sqrt \nu\, {\mathbb T}(\bar{v})$ and $\sqrt\nu \, {\mathbb T} (w)$ in $L^2(0,T;L^2(\Omega))$ 
and $\sqrt {\mu(\rho)}\in L^\infty(0,T;L^2(\Omega))$ to define in the usual way the first order derivative quantities. 

\proof
Thanks to \eqref{energyNSK} , we have
\begin{eqnarray}
\mathcal{E}_{NSK}(t)- \mathcal{E}_{NSK}(0) &   \leq  &  
 \displaystyle \int_{\Omega} \left(\dfrac{\rho}{2} |\bar{V}|^2-\rho \bar{v} \cdot \bar{V} + \dfrac{\rho}{2} |W|^2-  \rho w \cdot W\right)(t)   \nonumber    \\ 
 & {}  & - \displaystyle  \int_{\Omega}    \left(\dfrac{\rho}{2} |\bar{V}|^2-\rho \bar{v} \cdot\bar{V} + \dfrac{\rho}{2} |W|^2-  \rho w \cdot W\right)(0)   \nonumber \\
 &  {}  &- \displaystyle  \int_{\Omega} \left( H(r)+H'(r)(\rho-r)\right)(t)+\int_{\Omega}  \left( H(r)+H'(r)(\rho-r)\right)(0)   \nonumber  \\
  &  {}  &- \nu \displaystyle   \int_0^t  \int_{\Omega} \mu'(\rho)H''(\rho) |\nabla \rho |^2  +  \nu \left(I_1^{NS}+I_2^{NS}\right),
  \end{eqnarray}
  where
 \begin{eqnarray}
 \label{I1NS}
 2 I_1^{NS}
 &=& \int_0^t \int_\Omega
   \frac{\lambda(\rho)}{\sqrt{\mu(\rho)}}  \left(  | \dive \bar{V}|^2
     + | \dive W|^2   \right), \\
&&   -  2 \int_0^t \int_\Omega \frac{\lambda(\rho)}{\sqrt{\mu(\rho)}}
    \left(  
      {\rm Tr} ({\mathbb T}(\bar{v})) \dive \bar{V} 
   +  {\rm Tr} ({\mathbb T}(w)) \dive W \right), \nonumber 
 \end{eqnarray}
\begin{eqnarray}
\label{I2NS}
I_2^{NS}
&=&  \int_0^t \int_\Omega \mu(\rho)   \, \displaystyle \left(  |\nabla \bar{V}|^2
       + |\nabla W|^2 \right) \\
 &-  &  2 \int_0^t \int_\Omega \sqrt{\mu(\rho)}   \, \displaystyle \left( 
     {\mathbb T}(\bar{v}) : \nabla \bar{V} 
  +  {\mathbb T}(w) : \nabla W\right)
\end{eqnarray}
Using  \eqref{QNSNF2}, 
\begin{eqnarray}
\nonumber
  \partial_t(\rho \, w \cdot W)& = & \partial_t(\rho  \, w)\cdot W+\rho \, w \cdot \partial_t W \\
\nonumber
=  &  &  \left< - \mbox{div}( \rho \,  w \otimes u ) - \nabla p(\rho) 
    + \nu \,  \mbox{div}(\sqrt{\mu(\rho)} {\mathbb T}(w)) +
 \sqrt{\eps^2 - \nu^2} \, \mbox{div}(\sqrt{\mu(\rho)} {\mathbb T}(\bar{v}))  ;  W \right>_{W^{-2,1}(\Omega)\times W^{2,\infty}(\Omega)}  \nonumber \\
 &  &+  <A_1 ; W >_{W^{-2,1}(\Omega)\times W^{2,\infty}(\Omega)}
  +\rho \,  w \cdot \partial_t W, \nonumber 
 \end{eqnarray}
 where $$A_1=  \dfrac{\nu}{2} \nabla (\frac{\lambda(\rho)}{\sqrt{\mu(\rho)}} 
   {\rm Tr} ({\mathbb T}(w)) 
    + \dfrac{\sqrt{\eps^2 - \nu^2}}{2} \nabla (\frac{\lambda(\rho)}{\sqrt{\mu(\rho)}} 
   {\rm Tr} ({\mathbb T}(\bar{v})).$$
Using  \eqref{QNSNF3},
\begin{eqnarray}
\nonumber
 & &  \partial_t(\rho \,  \bar{v} \, \cdot \bar{V}) = \partial_t(\rho \,  \bar{v}) \,\cdot \bar{V}+\rho \,  \bar{v} \, \cdot \partial_t \bar{V} \\
  \nonumber 
 & = &   \left<-\mbox{div}(\rho \,  \bar v  \otimes u) + \nu \,  \mbox{div} (\sqrt{\mu(\rho)}
    {\mathbb T}(\bar{v}))  - \sqrt{\eps^2 - \nu^2} \, \mbox{div}( \sqrt{\mu(\rho)} ({\mathbb T}(w))^t) ;  \bar{V} 
 \right>_{W^{-2,1}(\Omega)\times W^{2,\infty}(\Omega)}     \\
 & & \hskip3cm + < A_2 ; \bar{V} >_{W^{-2,1}(\Omega)\times W^{2,\infty}(\Omega)}
    +\rho \,  \bar{v} \, \cdot \partial_t \bar{V}, 
\end{eqnarray}
where 
$$A_2=  \dfrac{\nu}{2} \,  \nabla (\frac{\lambda(\rho)}{\sqrt{\mu(\rho)}}
     {\rm Tr}({\mathbb T}(\bar{v})) )
- \dfrac{1}{2} \,  \sqrt{\eps^2 - \nu^2} \,   \nabla (\frac{\lambda(\rho)}{\sqrt{\mu(\rho)}}
     {\rm Tr}({\mathbb T}(w))).$$
Then,
Using \eqref{QNSNF1} 
\begin{eqnarray}
\nonumber
\int_0^t \int_{\Omega} \partial_t(H(r)+H'(r)(\rho-r))  & = &
\int_0^t \int_{\Omega}\left(  H'(r)  \,  \partial_t r +  \partial_t (H'(r)) (\rho-r)+H'(r) \,  \partial_t \rho  -H'(r) \,  \partial_t r \right) \\
\nonumber
& = & \int_0^t \int_{\Omega}\left(   \partial_t (H'(r)) (\rho-r)- H'(r) \mbox{div}\left(\rho \, u\right)  \right) \\
\nonumber
&  =  & \int_0^t \int_{\Omega}  \left( \partial_t (H'(r)) (\rho-r)+ \rho \nabla (H'(r))\cdot u  \right).
\end{eqnarray}
Since
$$\partial_t \left( \dfrac{\rho}{2} |\bar{V}|^2 \right)=\dfrac{1}{2} \partial_t \rho |\bar{V}|^2+ \rho \bar{V} \cdot \partial_t \bar{V}, \,  \quad
 \partial_t \left( \dfrac{\rho}{2} |W|^2 \right)=\dfrac{1}{2} \partial_t \rho |W|^2+ \rho W \cdot\partial_t W,$$
and since $\nabla \bar{v}, \nabla \bar{V}$ are symmetric matrices (recall that $v$ and $V$ are
gradient of functions), thanks to \eqref{QNSNF1} and integrations by parts we obtain
\begin{eqnarray}
\mathcal{E}_{NSK}(t)- \mathcal{E}_{NSK}(0)
  &  \leq  & \displaystyle \int_0^t \int_{\Omega} \rho \,  \partial_t \bar{V}  \cdot  (\bar{V} -  \bar v)
+ \displaystyle \int_0^t \int_{\Omega} \rho \,  (\nabla \bar{V} \, u) \cdot   (\bar{V} -  \bar v)   \nonumber  \\
& {}  &+ \int_0^t \int_{\Omega} \rho \, \partial_t W  \cdot   (W-w) 
+ \displaystyle \int_0^t \int_{\Omega} \rho \, ( \nabla W \, u)\cdot   (W-w) \nonumber  \\
 & {}  & + \nu \, \displaystyle  \int_0^t \int_{\Omega}
   \sqrt{\mu( \rho)} \,  \left( {\mathbb T}(\bar{v}) : \nabla \bar{V} +  {\mathbb T}(w)  : \nabla W \right)
\nonumber  \\
 & {}  &
 +  \sqrt{\eps^2 - \nu^2}  \,  \displaystyle \int_0^t \int_{\Omega}
 \sqrt{\mu( \rho )}\, \left(  {\mathbb T}(\bar{v})  : \nabla W   -  ({\mathbb T}(w))^t :  \nabla \bar{V}  \right)
 \nonumber  \\
& {}  &  - \displaystyle  \int_0^t \int_{\Omega}  p(\rho) \,  \mbox{div}W -  \nu \displaystyle   \int_0^t  \int_{\Omega} \mu'(\rho)H''(\rho) |\nabla \rho |^2
    \nonumber  \\
  &  {}  &  - \displaystyle \int_0^t \int_{\Omega} \partial_t(H'(r))(\rho-r)
-  \displaystyle \int_0^t \int_{\Omega} \rho\nabla( H'(r)) \cdot u   \nonumber  \\
  &  {}  &  +  \int_0^t \int_{\Omega} (\nu \, A_3 -\sqrt{\eps^2-\nu^2} \, A_4 ) +\nu  \left(I_1^{NS}+I_2^{NS}\right),   \nonumber  
 \end{eqnarray}
where
 $$2 \, A_3=    \frac{\lambda(\rho)}{\sqrt{\mu(\rho)}}  \,
       \left({\rm Tr}({\mathbb T}(\bar{v})) \dive \bar{V}  
           + {\rm Tr}({\mathbb T}(w))     \dive W  \right),$$
 $$2 \, A_4=     \frac{\lambda(\rho)}{\sqrt{\mu(\rho)}}  \, 
  \left( {\rm Tr}({\mathbb T}(w))   \dive \bar{V} 
   - {\rm Tr}({\mathbb T}(\bar{v})) \dive W   \right)$$
 and $I_1^{NS}$ and $I_2^{NS}$ are given through \eqref{I1NS}--\eqref{I2NS} which gives the proposition.
\cqfd

Let us introduce  that there exists a strong solution of 
\begin{eqnarray}
&\partial_t  r +\dive \left(r  \, U \right)=0,& \label{QNSentrel1}\\
\nonumber \\
&r \, \left(  \partial_t W  +   \nabla W \,U \right) +\nabla  p(r) - \nu \,  \dive (\mu(r)\,  \nabla W) -\dfrac{\nu}{2} \nabla(\lambda(r) \dive W)  \nonumber \\
 & =  \sqrt{\eps^2-\nu^2} \, \left(  \dive (\mu( r) \,   \nabla \bar{V})  +\dfrac{1}{2} \nabla(\lambda(r) \dive \bar{V})  \right), &\label{QNSentrel2} \\
\nonumber  \\
&r \, \left(  \partial_t \bar{V} +  \nabla   \bar{V} \, U  \right) - \nu \,  \dive (\mu(r) \, \nabla \bar{V}) -\dfrac{\nu}{2} \nabla(  \lambda(r) \dive \bar{V} )\nonumber \\
 & +\sqrt{ \eps^2 - \nu^2} \, \left(     \dive (\mu( r)\, ^t \nabla W)  +\dfrac{1}{2} \nabla(\lambda (r) \dive W) \right)
=0,& \label{QNSentrel3}
\end{eqnarray}
with $$U=W-\nu \, V, \quad \bar{V}=\sqrt{\eps^2-\nu^2} \, V$$
belonging to the class
\begin{eqnarray}
&&\nonumber 0 < {\rm inf}_{(0,T)\times \Omega}r  \le r \le {\rm sup}_{(0,T)\times \Omega} r < + \infty
\\
&& \nonumber   \nabla r \in L^2 (0,T;  L^\infty (\Omega)
         \cap L^1(0,T; W^{1, \infty}(\Omega))
 \\
&& \label{hypothesis}  W   \in    L^\infty (0,T; W^{2,\infty}(\Omega))  \cap 
                      W^{1,\infty}(0,T; L^\infty(\Omega)),
 \\
&& \nonumber \bar{V}    \in  L^\infty(0,T; W^{2,\infty}(\Omega))
                            \cap W^{1,\infty}(0,T; L^\infty(\Omega)),
 \\
&&\nonumber  \partial_t H'(r) \in L^1 (0,T;L^{\gamma/(\gamma-1)}(\Omega)), \qquad
       \nabla H'(r) \in L^1 (0,T; L^{2\gamma/(\gamma-1)}(\Omega))
\end{eqnarray}
and where $\bar{V}\vert_{t=0} = \sqrt{\varepsilon^2-\nu^2} \nabla \mu(r_0)/r_0$.
Defining
 \begin{eqnarray}
  I_3^{NS}= & {} &  \displaystyle   \int_0^t  \int_{\Omega}
   \lambda(\rho) \left(  (\dive \bar {V})^2 +(\dive W)^2 \right)   \label{I3NS}\\
   & & \nonumber  -   \int_0^t  \int_{\Omega}  \frac{\lambda(\rho)}{\sqrt{\mu(\rho)}}
    \left( {\rm Tr}({\mathbb T}({\bar v})) \dive \bar{V} 
         + {\rm Tr}({\mathbb T}(w))  \dive W \right) \\
   & + & \displaystyle   \int_0^t  \int_{\Omega}  \dfrac{\rho}{r} \Big( \nabla(\lambda(r) \dive \bar{V})\cdot(\bar V - \bar v )+ \nabla (\lambda(r) \dive W)\cdot(W-w) \Big), \nonumber  
 \end{eqnarray}
 \begin{eqnarray}
 I_4^{NS}= & {} &\displaystyle   \int_0^t  \int_{\Omega}  \mu(\rho) \left(|\nabla \bar{V}|^2+|\nabla W|^2 \right)  \label{I4NS}  \\
 & & \nonumber -   \int_0^t  \int_{\Omega}   \sqrt{\mu(\rho)} 
   \left({\mathbb T}({\bar v}) :  \nabla \bar{V} + \nabla W : {\mathbb T}(w) \right)  \\
 &  +  & \displaystyle   \int_0^t  \int_{\Omega}  \dfrac{\rho}{r} \Big(\dive \left( \mu(r)  \nabla \bar {V }\right)\cdot (\bar{V}-\bar{v}) + \dive\left(\mu(r) \nabla W\right)\cdot (W-w) \Big), \nonumber
 \end{eqnarray}
 \begin{eqnarray}
I_5^{NS} & = &
\int_0^t  \int_{\Omega}  \dfrac{\rho}{r} \Big(\dive \left( \mu(r)  \nabla \bar {V }\right)\cdot (W-w) - \dive\left(\mu(r) {}^t\nabla W\right)\cdot (\bar{V}-  \bar{v}) \Big) \label{I5NS} \\
&  & \nonumber +   \int_0^t\int_{\Omega}
  \sqrt{\mu( \rho )} \left({\mathbb T}({\bar v}) : \nabla W -   ({\mathbb T}(w))^t :\nabla \bar{V} \right),
 \end{eqnarray}
 \begin{eqnarray}
2 \, I_6^{NS}  & = &
  \int_0^t \int_{\Omega}  \dfrac{\rho}{r} \left( \nabla (\lambda(r) \dive \bar{V})\cdot (W-w)  - \nabla (\lambda(r) \dive W) \cdot (\bar{V}-\bar{v}) \right)  \label{I6NS} \\
 & &- \int_0^t \int_{\Omega} \frac{\lambda(\rho)}{\sqrt{\mu(\rho)}} \,  
  \left({\rm Tr}({\mathbb T}(w)) \dive \bar{V} - 
  {\rm Tr}( {\mathbb T}({\bar v})) \dive W   \right),  \nonumber 
  \end{eqnarray}
we have
 \begin{prop} \label{prop20}
 Let $(r,\bar{V},W)$ be  a  strong solution  of \eqref{QNSentrel1}-\eqref{QNSentrel3}
belonging to the class \eqref{hypothesis}.  Let us assume that $\bar{V_0} = 
\sqrt{\varepsilon^2 - \nu^2} \nabla \mu(r_0)/r_0$.
Any weak  solution $(\rho,\bar{v},w)$ of the augmented system \eqref{QNSNF1}--\eqref{QNSNF-constraint} satisfies the following  inequality
  \begin{eqnarray*}
\mathcal{E}_{NSK}(t)- \mathcal{E}_{NSK}(0) 
&  \leq &
 \displaystyle   \int_0^t \int_{\Omega}  \rho \, \left[   ( \nabla  \bar{V} \, (u-U)) \cdot  (\bar{V}-\bar{v}) 
 +  (\nabla W \, (u-U) ) \cdot  (W-w)  \right]  \nonumber  \\
 & &  -  \displaystyle   \int_0^t \int_{\Omega} \left( p(\rho) -p(r) - p'(r) \,   (\rho-r)  \right) \dive  U
 \nonumber  \\
&  & +   \displaystyle\dfrac{\nu}{\sqrt{\varepsilon^2-\nu^2}}    \int_0^t \int_{\Omega} \left(\rho \, \nabla (H'(r)) \cdot (\bar{v}-\bar{V})- p(\rho) \, \dive \bar{V} \right) \\
 & & - \nu \displaystyle   \int_0^t  \int_{\Omega} \mu'(\rho)H''(\rho) |\nabla \rho |^2 \nonumber  \\
  &  & + \dfrac{\nu}{2} I_3^{NS} +\displaystyle  \nu I_4^{NS} +  \sqrt{\eps^2-\nu^2} \, \left(  I_5^{NS}+ I_6^{NS} \right),  \nonumber \\
\end{eqnarray*}
where $I_{i}^{NS}$ for $i=3,4,5,6$ are given by \eqref{I3NS}-\eqref{I6NS}.
\end{prop}

\proof
Multiplying \eqref{QNSentrel2} by  $\dfrac{\rho}{r}(W-w)$ and \eqref{QNSentrel3} by $\dfrac{\rho}{r}(\bar{V}-\bar{v})$, integrating with respect to time and space, and using \eqref{QNSentrel1}, we obtain
  \begin{eqnarray*}
   \mathcal{E}_{NSK}(t)- \mathcal{E}_{NSK}(0) & \leq & \nu   \int_0^t \int_{\Omega}  \dfrac{\rho}{r} \,  \left[ \mbox{div}(\mu( r ) \nabla \bar{V})\cdot (\bar{V}-\bar{v})
+\dfrac{1}{2}  \nabla (\lambda(r)\dive \bar{V}) \cdot (\bar{V}-\bar{v}) \right]
 \nonumber \\ 
 &  & -  \sqrt{\eps^2-\nu^2} \int_0^t \int_{\Omega} \dfrac{\rho}{r} \,  \left[   \mbox{div}(\mu( r ){}^t\nabla W) \cdot (\bar{V}-\bar{v}) 
 +\dfrac{1}{2}  \nabla(\lambda(r) \dive W) \cdot (\bar{V}-\bar{v}) \right]   \nonumber \\ 
 &   &+  \int_0^t\int_{\Omega}  \rho  \,(\nabla \bar{V}\, (u-U)) \cdot  (\bar{V}-\bar{v})
 +  \int_0^t \int_{\Omega}  \rho  \,( \nabla W \, (u-U)) \cdot  (W-w)
\nonumber \\
&  & +  \nu  \int_0^t\int_{\Omega}  \dfrac{\rho}{r}  \, \left[ \mbox{div}(\mu( r )\nabla W)\cdot (W-w)
 +\dfrac{1}{2} \, \nabla (\lambda(r) \dive W)\cdot (W-w)\right]  
 \end{eqnarray*}
 \begin{eqnarray*}
&  &
+  \sqrt{\eps^2-\nu^2}\int_0^t \int_{\Omega} \dfrac{\rho}{r}  \left[ \mbox{div}(\mu( r ) \nabla \bar{V})\cdot (W-w)
  +\dfrac{1}{2}  \nabla( \lambda(r) \dive \bar{V})\cdot(W-w)\right] \nonumber \\
  & 	  & +   \nu \displaystyle   \int_0^t  \int_{\Omega}\mu(\rho) \,  
       \left( |\nabla \bar{V}|^2 +|\nabla W|^2\right)  \\
   & &  -     \nu   \int_0^t  \int_{\Omega} 
       \sqrt{\mu(\rho)} \left( {\mathbb T}(\bar{v}) : \nabla \bar{V} 
       +  {\mathbb T}(w): \nabla W   \right) 
         \nonumber \\
 &  & +   \sqrt{\eps^2 - \nu^2} \displaystyle   \int_0^t \int_{\Omega}
 \sqrt{\mu(\rho)} \left( {\mathbb T}(\bar{v}):  \nabla W -
      ( {\mathbb T}(w))^t: \nabla \bar{V}
 \right)   \nonumber  \\
 \end{eqnarray*}
 \begin{eqnarray*}
  & &  +  \sqrt{\eps^2 - \nu^2}  \displaystyle   \int_0^t \int_{\Omega}
\dfrac{ \lambda(\rho)}{2\sqrt{\mu(\rho)}} 
  \left({\rm Tr}({\mathbb T}(\bar{v})) \dive W - 
  {\rm Tr}({\mathbb T}(w))  \dive \bar{V}  \right)  \nonumber  \\
   &  & -  \dfrac{\nu}{2}  \displaystyle   \int_0^t \int_{\Omega} \lambda(\rho) \left( (
    \dive \bar{V})^2  + ( \dive W )^2  \right)  \nonumber  \\
     &  & +  \dfrac{\nu}{2}  \displaystyle   \int_0^t \int_{\Omega} 
     \frac{\lambda(\rho)}{\sqrt{\mu(\rho)} }\left(  
      {\rm Tr}({\mathbb T}(\bar{v})) \dive\bar{V} 
       + {\rm Tr}({\mathbb T}(w)  \dive W \right)  \nonumber  \\
  &   &  +  \displaystyle   \int_0^t  \int_{\Omega}  \left( p(r) \dive U -p(\rho) \dive W  \right)-  \nu \displaystyle   \int_0^t  \int_{\Omega}\mu'(\rho) H''(\rho) |\nabla \rho |^2 +I_7^{NS},   \nonumber   \\
\end{eqnarray*}
where
\begin{eqnarray}
I_7^{NS}& = &  -\int_0^t \int_\Omega \dfrac{\rho}{r}  \nabla (p(r)) \cdot (W-w)-\int_0^t \int_{\Omega} \partial_t (H'(r))(\rho-r)  \label{I7NS}\\
 & & -\int_0^t\int_{\Omega} \rho\nabla (H'(r)) \cdot u + \int_0^t\int_{\Omega}H'(r) \, \partial_t r. \nonumber
 \end{eqnarray}
Using \eqref{QNSentrel1} $H'(r)\partial_t r+ H'(r)\mbox{div}(r \, U)=0$ which leads, with an integration by parts, to
$$\displaystyle \int_0^t \int_\Omega ( H'(r)\partial_t r - r \nabla (H'(r)) \cdot U)  =0.$$
Then 
$$\displaystyle \int_0^t \int_\Omega(  H'(r)\partial_t r -   \nabla p(r) \cdot U ) =0,$$
or 
\begin{equation}
\label{proofprop15_1}
\int_0^t \int_\Omega ( H'(r)\partial_t r + p(r) \mbox{div}(U)) =0.
\end{equation}
\noindent
Moreover
$$\partial_t (H'(r))=-p'(r) \mbox{div}U-H''(r )\nabla r \cdot U=-p'(r) \mbox{div}U -\nabla (H'(r)) \cdot U.$$
Then
\begin{eqnarray*}
I_7^{NS} & = &  \int_0^t \int_\Omega  \rho \, \nabla (H'(r) ) \cdot (- W+w +U -u)+ \int_0^t \int_{\Omega} p'(r) \mbox{div}U (\rho-r)  \\
& = & \nu  \int_0^t \int_\Omega \rho  \nabla( H'(r) ) \cdot (v-  V)+ \int_0^t \int_{\Omega} p'(r) \mbox{div}U (\rho-r). \\
 \end{eqnarray*}
Therefore 
\begin{eqnarray}
\mathcal{E}_{NSK}(t)- \mathcal{E}_{NSK}(0) 
&  \leq &
 \displaystyle   \int_0^t \int_{\Omega} \rho \,  \left[  ( \nabla  \bar{V} \, (u-U)) \cdot  (\bar{V}-\bar{v})
 + ( \nabla W  \, (u-U)) \cdot  (W-w)  \right]  \nonumber  \\
  &  & +   \displaystyle   \int_0^t \int_{\Omega} \left[p(r) \,  \mbox{div} U-p(\rho) \, \mbox{div}U -\nu \,  p(\rho) \, \mbox{div}V  \right]
 - \nu \displaystyle   \int_0^t  \int_{\Omega}\mu'(\rho)H''(\rho) |\nabla \rho |^2 
 \nonumber  \\
 &  &  +  \displaystyle   \dfrac{\nu}{2}   \int_0^t\int_{\Omega}  \lambda(\rho)  \left( (\dive \bar{V})^2 + (\dive W)^2  \right) \nonumber \\
  &  &  -  \displaystyle   \dfrac{\nu}{2}   \int_0^t\int_{\Omega}
     \frac{\lambda(\rho)}{\sqrt{\mu(\rho)} } \left( 
  {\rm Tr}({\mathbb T}(\bar{v})) \dive \bar{V}  +
  {\rm Tr}({\mathbb T}(w))   \dive W \right) \nonumber \\
 &  & +   I_7^{NS}+  \nu \,  I_4^{NS}+ \nu I_8^{NS}+ \sqrt{\eps^2-\nu^2} \,  ( I_5^{NS} +  I_6^{NS})  \nonumber 
\end{eqnarray}
where 
\begin{equation}
\label{I8NS}
I_8^{NS}=\dfrac{1}{2}\int_0^t \int_\Omega\dfrac{\rho}{r} \left(\nabla(\lambda(r)\dive \bar{V}) \cdot (\bar{V}-\bar{v})+ \nabla(\lambda(r)\dive W) \cdot (W-w)\right),
\end{equation}
and $I_{i}^{NS}$ for $i=4, 5, 6, 7$ are given by \eqref{I4NS}-\eqref{I7NS}.
Finally
\begin{eqnarray*}
\mathcal{E}_{NSK}(t)- \mathcal{E}_{NSK}(0) 
&  \leq &
 \displaystyle   \int_0^t \int_{\Omega} \left[  \rho \, ( \nabla  \bar{V} \, (u-U)) \cdot  (\bar{V}-\bar{v})
 +\rho \,(  \nabla W \, (u-U)  ) \cdot  (W-w)\right]  \nonumber  \\
 &  &  +  \displaystyle   \int_0^t \int_{\Omega} \left[(p'(r)  (\rho-r)-p(\rho) +p(r) )\,  \dive  U-\nu \,  p(\rho) \, \dive V \right]
 \nonumber  \\
&  &  - \dfrac{\nu}{\sqrt{\eps^2-\nu^2}}  \displaystyle     \int_0^t \int_{\Omega} \rho\nabla( H'(r)) \cdot (\bar{V}-\bar{v})
- \nu \displaystyle   \int_0^t  \int_{\Omega}\mu'(\rho) H''(\rho) |\nabla \rho |^2 \nonumber  \\
  &  & +  \dfrac{\nu }{2}  \displaystyle   \int_0^t  \int_{\Omega} \lambda(\rho) \left(  (\dive \bar {V})^2  +(\dive W)^2  \right) \nonumber  \\ 
  &  & -  \dfrac{\nu }{2}  \displaystyle   \int_0^t  \int_{\Omega}
\frac{\lambda(\rho)}{\sqrt{\mu(\rho)}}
  \left(  
   {\rm Tr}({\mathbb T}(\bar{v}))\dive \bar{V} 
   +  {\rm Tr}({\mathbb T}(w)) \dive W \right) \nonumber  \\ 
   &  & +  \dfrac{\nu}{2} \displaystyle   \int_0^t  \int_{\Omega}  \dfrac{\rho}{r} \left( \nabla(\lambda(r) \dive \bar{V})\cdot(\bar V - \bar v )+ \nabla (\lambda(r) \dive W)\cdot(W-w )\right) \nonumber  \\
  &  & +  \nu \displaystyle   \int_0^t  \int_{\Omega}  \mu(\rho) \left( |\nabla \bar{V}|^2+|\nabla W|^2 \right)  \nonumber  \\
   &  & -  \nu \displaystyle   \int_0^t  \int_{\Omega} 
    \sqrt{\mu(\rho)} \left( 
      {\mathbb T}({\bar v}) : \nabla \bar{V}  +\nabla W : {\mathbb T}(w) \right)  \nonumber  \\
&   & + \nu  \displaystyle   \int_0^t  \int_{\Omega}  \dfrac{\rho}{r} \Big(\dive\left( \mu(r)  \nabla \bar {V }\right)\cdot(\bar{V}-\bar{v}) + \dive\left(\mu(r) \nabla W\right)\cdot (W-w)\Big)\\
  &  & +   \displaystyle   \sqrt{\eps^2-\nu^2} \, \left(  I_5^{NS}+ I_6^{NS} \right),  \nonumber \\
\end{eqnarray*}
with $I_5^{NS}$ and $I_6^{NS}$ given by \eqref{I5NS} and \eqref{I6NS}. This gives the proposition.
\cqfd

\begin{lem}\label{lem_I5I6}
Let $I_5^{NS}$ given by \eqref{I5NS} and $I_6^{NS}$ given by \eqref{I6NS}. Under the assumptions of Proposition \ref{prop20},  we have
 \begin{eqnarray*}
I_5^{NS} & = &
- \int_0^t\int_{\Omega} \rho \left( \dfrac{\mu(\rho)}{\rho}-\dfrac{\mu(r)}{r} \right) \Big(\dive(\nabla \bar V)\cdot (W-w)+\dive({}^t\nabla W)\cdot (\bar{v}-\bar V)\Big) \\
  & & -  \int_0^t\int_{\Omega} \rho \, \Big(\nabla \bar V (W-w) +\nabla W(\bar v -\bar V)\Big) \cdot (v-V),
 \end{eqnarray*}
 and
 \begin{eqnarray*}
2 \,  I_6^{NS} & =  & -2 \,  \int_0^t  \int_{\Omega}  \rho \left(   \mu''(\rho) \nabla \rho  -  \mu''(r)\nabla r \right) \cdot \left(  (W-w) \, \dive \bar{V} + (\bar{v}-\bar{V}) \, \dive W \right)  \\
&  & -  \int_0^t  \int_{\Omega}  \left( \lambda (\rho) -\dfrac{\rho}{r}\lambda (r)\right)  \left(  (W-w) \, \cdot \nabla( \dive \bar{V}) + (\bar{v}-\bar{V}) \,\cdot \nabla ( \dive W )\right).   \\
  \end{eqnarray*}
\end{lem}

\proof The proof follows the same lines that the ones for \eqref{I1} in the Euler-Korteweg
section.
\cqfd

\begin{lem}\label{lem22}
Let $I_5^{NS}$ given by \eqref{I5NS} and $I_6^{NS}$ given by \eqref{I6NS}. Let us assume  $\mu(\rho)=\rho^{(s+3)/2}$ with $ \gamma \geq s+2$ and $s \geq -1$. Under the assumptions of Proposition \ref{prop20} we have
\begin{eqnarray*}
\Bigl| I_5^{NS}+I_6^{NS}  \Bigr|
 & \leq & 
 C  \,   \int_0^t  \int_{\Omega}  \left( \rho   |w -W |^2 +    \rho   |v -V |^2 +  \rho   |\bar{v} -\bar{V} |^2  + H(\rho | r)  \right),
\end{eqnarray*}
where $C=C(r, \bar V, W)$ is  a uniformly bounded constant on $\R^+ \times \Omega$.
\end{lem}

\proof
By definition of $\lambda(\rho)$, lemma \ref{lem_I5I6} directly leads to
  \begin{eqnarray*}
I_5^{NS}+I_6^{NS} & = & 
- \int_0^t  \int_{\Omega} \rho \, (  \nabla  \bar{V}(  W-w )) \cdot (v-V)   +  \int_0^t  \int_{\Omega} \rho  \,( \nabla W (\bar{V}-\bar{v})  )\cdot (v-V)  \\
   & &  -  \int_0^t  \int_{\Omega}  \rho \left(   \mu''(\rho) \nabla \rho  -  \mu''(r)\nabla(r) \right) \cdot \left(  (W-w) \, \dive \bar{V} + (\bar{v}-\bar{V}) \, \dive W  \right)  \\
& &  - \int_0^t  \int_{\Omega} \rho  (  \mu'(\rho)  -  \mu'(r)) \left(  (W-w) \cdot \nabla ( \dive \bar{V} ) +(\bar{v}-\bar{V}) \cdot \nabla ( \dive W )\right).
 \end{eqnarray*}
Moreover in an analogous way than for lemma \ref{lem6}, we can show that
$$\left| \int_0^t  \int_{\Omega}  \rho \left(   \mu''(\rho) \nabla \rho  -  \mu''(r)\nabla(r) \right) \cdot \left(  (W-w) \, \dive \bar{V} + (\bar{v}-\bar{V}) \, \dive W  \right) \right|$$
$$\leq  C  \,    \int_0^t  \int_{\Omega} \rho \left(|W-w|^2 + |\bar v-\bar V|^2  \right)$$
and
$$\left|\int_0^t  \int_{\Omega} \rho  (  \mu'(\rho)  -  \mu'(r)) \left(  (W-w) \cdot \nabla ( \dive \bar{V} ) +(\bar{v}-\bar{V}) \cdot \nabla ( \dive W )\right)\right|$$
$$\leq  C  \int_0^t   \int_\Omega  \left( H(\rho | r) +  \rho (|w-W|^2 +|\bar v-\bar V|^2) \right).$$
\cqfd

\begin{lem}
Let $I_3^{NS}$ given by \eqref{I3NS} and $I_4^{NS}$ given by \eqref{I4NS}. Under the assumptions of Proposition \ref{prop20}, we have
\begin{eqnarray*}
I_3^{NS} & = &  -2 \int_0^t  \int_{\Omega} \, \rho( \mu''(\rho)\nabla \rho  - \mu''(r) \nabla r) \cdot (( \bar V- \bar v)\dive \bar V + (W-w)\dive W ) \\
 & & -\int_0^t  \int_{\Omega}\left(\lambda(\rho)-\dfrac{\rho}{r} \lambda(r)\right)(\nabla (\dive \bar V)\cdot( \bar V-\bar v)+\nabla (\dive W)\cdot (W-w)),
\end{eqnarray*}
and 
\begin{eqnarray*}
I_4^{NS} & = &  - \int_0^t  \int_{\Omega}\rho ( v -V) \cdot \left(\nabla \bar V (\bar V-\bar v)+ \,{}^t\nabla W (W-w)\right)  \\
 & & -\int_0^t  \int_{\Omega} \rho\left(\dfrac{\mu(\rho)}{\rho} -\dfrac{ \mu(r)}{r} \right)(\dive (\nabla \bar V)\cdot(\bar V-\bar v)+\dive ( \nabla W)\cdot(W-w)).
\end{eqnarray*}
\end{lem}

\proof The proof follows the same lines that the ones for \eqref{I1} in the Euler-Korteweg
section.
\cqfd

Using the previous lemma and the symmetry of $\nabla \bar{V}$, we obtain the following lemma
\begin{lem}\label{lem26}
Let $I_3^{NS}$ given by \eqref{I3NS} and $I_4^{NS}$ given by \eqref{I4NS}. We assume $\mu(\rho)=\rho^{(s+3)/2}$ with $ \gamma \geq s+2$ and $s \geq -1$. Under the assumptions of Proposition \ref{prop20}, we have
\begin{eqnarray*}
\Bigl| \dfrac{1}{2}I_3^{NS}+I_4^{NS}\Bigr| & \leq & 
  C   \int_0^t  \int_{\Omega} \left(  \rho   |v -V|^2 +  \rho   |\bar{v} -\bar{V}|^2 
+  \rho   |w -W|^2 + H(\rho|r)  \right),  
\end{eqnarray*}
where $C=C(r, \bar V, W)$ is  a uniformly bounded constant on $\R^+ \times \Omega$.
\end{lem}

\proof
We have:
\begin{eqnarray*}
\dfrac{1}{2}I_3^{NS}+I_4^{NS} & = & 
 - \int_0^t  \int_{\Omega}\rho \,  ( v -V) \cdot \left(\nabla \bar V (\bar V-\bar v)+ \,{}^t\nabla W (W-w)\right)  \\
& & - \int_0^t  \int_{\Omega} \, \rho \, ( \mu''(\rho)\nabla \rho  - \mu''(r) \nabla r) \cdot (( \bar V- \bar v)\dive \bar V + (W-w)\dive W ) \\
& & - \int_0^t  \int_{\Omega} \, \rho \, ( \mu'(\rho)  - \mu'(r) )\left( \nabla (\dive \bar V)\cdot(\bar V-\bar v)+\nabla( \dive W)\cdot (W-w)\right) \\
& &-\int_0^t  \int_{\Omega} \rho \,  \left(\dfrac{\mu(\rho)}{\rho} -\dfrac{\mu(r)}{r}\right)\left(\dive(\nabla W)-\nabla(\dive W)\right)\cdot (W-w).
\end{eqnarray*}
In an analogous way than for the lemma \ref{lem6}, we can show
 $$\int_0^t  \int_{\Omega} \rho \left(\dfrac{\mu(\rho)}{\rho} -\dfrac{\mu(r)}{r}\right)\left(\dive(\nabla W)-\nabla(\dive W)\right)\cdot (W-w) 
\leq   C  \int_0^t\int_\Omega \left( H(\rho | r)  +  \rho \,  |W-w|^2 \right).$$
Then using an analogous result than the one used in the proof of lemma \ref{lem22} we obtain the result.
\cqfd

Let us now define
\begin{eqnarray*}
 I_{11}^{NS} & =  & -  \dfrac{\nu}{\sqrt{\eps^2-\nu^2}}\int_0^t \int_{\Omega} \left(  \rho  \nabla( H'(r)) \cdot (\bar{V}-\bar{v}) + p(\rho) \dive V \right)  - \nu \int_0^t \int_{\Omega} \mu'(\rho) H''(\rho) |\nabla \rho|^2. 
 \end{eqnarray*}
 Using the definition of $H$ and an integration by parts, we obtain
 \begin{eqnarray}
 I_{11}^{NS} & = & - \dfrac{\nu}{\eps^2-\nu^2}\int_0^t \int_{\Omega} \rho \,  \left(  \dfrac{p'(r)}{\mu'(r)} \bar{V}  -  \dfrac{p'(\rho)}{\mu'(\rho)} \bar{v}\right) \cdot ( \bar{V}-\bar{ v }), \label{I11NS}
 \end{eqnarray}
 with $v=\nabla(\mu(\rho))/\rho,~\bar{v}=\sqrt{\eps^2-\nu^2}v,~V=\nabla(\mu(r))/r,~\bar{V}=\sqrt{\eps^2-\nu^2}V$.
We can show the following proposition.
\begin{prop}\label{prop28}
Let $I_{11}^{NS}$ given by \eqref{I11NS}. Assuming  $\mu(\rho)=\rho^{(s+3)/2} $ with $ \gamma  \geq s+2$, $s \geq -1$ and the hypothesis of Proposition \ref{prop20},
there exists a contant $C=C(r,U,\bar{V},W)$ uniformly bounded on $\R^+ \times \Omega$ such that
\begin{eqnarray*}
I_{11}^{NS}  \le  
 C  \dfrac{\nu}{\eps^2-\nu^2}\int_0^t \int_{\Omega}  H(\rho|r).
 \end{eqnarray*}
\end{prop}

\proof
Using Lemma \ref{lem_tech5}, we can write 
\begin{equation}\label{I11}
I_{11}^{NS} = - 
  \dfrac{\nu}{\eps^2-\nu^2} \int_0^t\int_{\Omega}  \rho  \, \dfrac{p'(\rho)}{\mu'(\rho)}|\bar{V}-\bar{v}|^2- I_{12}^{NS}
  \end{equation}
where
\begin{equation}\label{I12}
I_{12}^{NS} = \dfrac{\nu}{\eps^2-\nu^2} \displaystyle \int_0^t\int_\Omega \left(\sqrt{\eps^2-\nu^2}\nabla \phi_1(\rho|r)+\phi_2(\rho|r)\bar{V} \right)\cdot \bar{V}.
\end{equation}
Using an integration by parts
\begin{eqnarray*}
I_{12}^{NS} &= &-\dfrac{\nu}{\sqrt{\eps^2-\nu^2}}\int_0^t \int_\Omega \phi_1(\rho|r)\dive( \bar{V})+ \dfrac{\nu}{\eps^2-\nu^2}\int_0^t \int_\Omega \phi_2(\rho|r)\bar{V} \cdot \bar{V}.
\end{eqnarray*}
Now using lemma \ref{lem_tech3} we obtain
\begin{eqnarray*}
 I_{12}^{NS} \le  \dfrac{C \, \nu}{\sqrt{\eps^2-\nu^2}} \int_0^t \int_\Omega H(\rho|r) + \dfrac{C \, \nu}{\eps^2-\nu^2} \int_0^t \int_\Omega H(\rho|r) \leq  \dfrac{C \, \nu}{\eps^2-\nu^2} \int_0^t \int_\Omega H(\rho|r),
 \end{eqnarray*}
which gives the result due to the expression \eqref{I11} and the sign of the first quantity
in the right-hand side.
\cqfd

\begin{theorem} \label{th17}
Assuming $\mu(\rho)=\rho^{(s+3)/2}$, $\gamma \geq s+2$ and $s \geq -1$, any weak  solution $(\rho,\bar{v},w)$ of System \eqref{QNSNF1}-\eqref{QNSNF-constraint} satisfies the following  inequality
 \begin{eqnarray} 
\mathcal{E}_{NSK}(t) -\mathcal{E}_{NSK}(0) &  \leq  &    C  \left( 1+ \dfrac{\nu}{\eps^2-\nu^2} \right)  \int_0^t   \mathcal{E}_{{NSK}} (\xi)d\xi
\end{eqnarray}
 where $(r,\bar{V},W)$ is a  strong solution of \eqref{QNSentrel1}-\eqref{QNSentrel3}
belonging to the class \eqref{hypothesis} and where $C=C(r,U,\bar{V},W)$  is a constant uniformly bounded on $\R^+ \times \Omega$.
\end{theorem}

\proof
Thanks to Proposition \ref{prop20} we have
\begin{eqnarray*}
\mathcal{E}_{NSK}(t)- \mathcal{E}_{NSK}(0) 
&  \leq &
 \displaystyle   \int_0^t \int_{\Omega}  \rho \, \left[   ( \nabla  \bar{V} \, (u-U)) \cdot  (\bar{V}-\bar{v}) 
 +  (\nabla W \, (u-U) ) \cdot  (W-w)  \right]  \nonumber  \\
 & &  -  \displaystyle   \int_0^t \int_{\Omega} \left( p(\rho) -p(r) - p'(r) \,   (\rho-r)  \right) \dive  U
 \nonumber  \\
&  &+  I_{11}^{NS}  + \dfrac{\nu}{2} I_3^{NS} +\displaystyle  \nu I_4^{NS} +  \sqrt{\eps^2-\nu^2} \, \left(  I_5^{NS}+ I_6^{NS} \right),  \nonumber 
\end{eqnarray*}
with $I_{i}^{NS}$ for $i=3, 4, 5, 6$ given by \eqref{I3NS}-\eqref{I6NS} and $I_{11}^{NS}$ given by \eqref{I11NS}. This gives with  the regularity of $U, \, \bar{V}$ and $W$ and the previous lemmas
\begin{eqnarray*}
\mathcal{E}_{NSK}(t)- \mathcal{E}_{NSK}(0) 
&  \leq &
 C  \displaystyle   \int_0^t \int_{\Omega}  \rho \ \left(  |u-U|^2 + |\bar{v}-\bar{V}|^2 + |w-W|^2  \right) 
  \nonumber  \\
 & &  -  \displaystyle   \int_0^t \int_{\Omega} \left( p(\rho) -p(r) - p'(r) \,   (\rho-r)  \right) \dive  U
 \nonumber  \\
&  &+   C \dfrac{\nu}{\eps^2-\nu^2}  \int_0^t \int_{\Omega}  H(\rho|r) \\
 &  \leq &  C \left( 1+ \dfrac{\nu}{\eps^2-\nu^2} \right)   \displaystyle  \int_0^t \mathcal{E}_{NSK}(\xi)d\xi.
\end{eqnarray*}
\cqfd

\begin{cor}
Let $(r,\bar{V},W)$ be a  strong solution  of \eqref{QNSentrel1}-\eqref{QNSentrel3}
in the class belonging to the class \eqref{hypothesis}. Assuming $\mu(\rho)=\rho^{(s+3)/2}$, $\gamma \geq s+2$ and $s \geq -1$ any weak  solution $(\rho,w,\bar{v})$ of \eqref{QNSNF1}-\eqref{QNSNF-constraint} satisfies the following  ine\-quality
\begin{eqnarray}\label{equa28}
\nonumber
\mathcal{E}_{NSK}(\rho,\bar{v},w |  r, \bar{V},W)(t)  \leq  &  \mathcal{E}_{NSK}(\rho,\bar{v},w |  r, \bar{V},W)(0)
\exp  \left(C\left(1+\dfrac{\nu}{ \eps^2-\nu^2} \right)\, t\right),
 \end{eqnarray}
where $C=C(r,U,\bar{V},W)$  is a constant uniformly bounded on $\R^+ \times \Omega$.
\end{cor}

\proof
Thanks to the previous proposition and the Gronwall's Lemma, we have the inequality.
\cqfd

Let $U$ be a given and smooth function. We define $r$ as the strong solution of \eqref{NSQF1}, and
 we introduce the functions $\mathscr{E}_1^\nu$ and $\mathscr{E}_2^\nu$ such that 
\begin{eqnarray}
&      \mathscr{E}_1^\nu(r,\bar{V},W),=   r \, \left(  \partial_t W  +  U \cdot \nabla W \right) +\nabla p(r) - \nu \,  \dive (\mu(r) \,  \nabla W) -\dfrac{\nu}{2} \nabla(\lambda(r) \dive W)   \label{43}   \\
 & -  \sqrt{\eps^2-\nu^2} \,  \left(  \dive (\mu(r) \,   \nabla \bar{V}) +\dfrac{1}{2} \nabla(\lambda(r) \dive \bar{V}) \right)  \nonumber  \\
& 0= r \, \left(  \partial_t \bar{V} + U  \cdot \nabla   \bar{V}  \right) - \nu \,  \dive (\mu(r) \nabla \bar{V}) -\dfrac{\nu}{2} \nabla(  \lambda(r) \dive \bar{V} ) \, \label{44} \\
& +\sqrt{ \eps^2 - \nu^2} \,   \left(   \dive (\mu(r)\, ^t \nabla W)  +\dfrac{1}{2} \nabla(\lambda (r) \dive W) \right)
 \nonumber  
\end{eqnarray}
where $\bar{V}= \sqrt{\varepsilon^2-\nu^2} \nabla \mu(r)/r$.
In a same way than for the proof of Theorem \ref{th17}, we have the following result.
  \begin{prop} \label{prop19} Let us assume $\mu(\rho)=\rho^{(s+3)/2}$ ({\it i.e.} $K(\rho)=\dfrac{(s+3)^2}{4}\rho^s$), $\gamma \geq s+2$ and $s \geq -1$.
Let $(\rho,\bar{v},w)$ be a global weak  solution of System \eqref{QNSNF1}-\eqref{QNSNF-constraint} and $(r,\bar{V},W)$ a strong solution of   \eqref{NSQF1}, \eqref{43}-\eqref{44}
in the clas \eqref{hypothesis}.
Then \begin{eqnarray} 
\mathcal{E}_{{NSK}}(\rho,\bar{v},w |  r, \bar{V},W)(t) - \mathcal{E}_{{NSK}}(\rho,\bar{v},w |  r, \bar{V},W)(0)
 &  \leq  &C \left( 1 +\dfrac{\nu}{\eps^2-\nu^2} \right)  \int_0^t \mathcal{E}_{{NSK}}
+b^\nu(t),\nonumber
\end{eqnarray}
with
$$b^\nu(t)=\int_0^t \int_{\Omega} \left[ \dfrac{\rho}{r} \mathscr{E}_1^\nu \cdot (W-w)\right],$$
and where $C=C(r,U,\bar{V},W)$  is a constant uniformly bounded on $\R^+ \times \Omega$.
 \end{prop}

Using the Gronwall's Lemma, we immediately obtain the following corollary.
\begin{cor}\label{corsoldisNSQ}
Let us assume $\mu(\rho)=\rho^{(s+3)/2}$ ({\it i.e.} $K(\rho)=\dfrac{(s+3)^2}{4}\rho^s$) with $\gamma \geq s+2$ and $s \geq -1$.
Let $(\rho,\bar{v},w)$ be a weak  solution of System \eqref{QNSNF1}-\eqref{QNSNF-constraint}  and $(r,\bar{V},W)$ a strong solution of   \eqref{NSQF1}, \eqref{43}-\eqref{44}
in the class \eqref{hypothesis}.
Then
\begin{eqnarray}
\mathcal{E}_{{NSK}}(t) 
 &  \leq  & \mathcal{E}_{{NSK}}(0) \, \exp(F^\nu \, t)+F^\nu  \int_0^t b^\nu(\xi) \, \exp(F^\nu \, (t-\xi)) \, d\xi
+b^\nu(t),\nonumber
\end{eqnarray}
where  $b^\nu$ is defined in Proposition \ref{prop19} and
$$F^\nu=C \left( 1 +\dfrac{\nu}{\eps^2-\nu^2} \right),$$ 
with $C=C(r,U,\bar{V},W)$ a constant uniformly bounded on $\R^+\times\Omega$.
\end{cor}


\subsection{Dissipative solution and weak-strong uniqueness result}

Let us now give the definition of what is called a dissipative solution of the compressible Navier-Stokes-Korteweg System.  To this end, let $U$ be a smooth function, then
$(r,{\mathcal E}^\nu(r,U))$ defined through Equation \eqref{NSQF1} and 
\begin{eqnarray}
\label{soldisQNS1}
\mathscr{E}^\nu(r,U) =
r \left( \partial_t U+ U \cdot  \nabla U  \right)+\nabla p(r)-2 \, \nu \,  \hbox{div}(\mu(r)D(U))
 -  \nu \nabla (\lambda(\rho)\dive U)   \nonumber   \\
 +
\eps^2\left[ ( \dive(\mu(r)^t \nabla V )+\dfrac{1}{2}  \nabla(\lambda(r) \dive V)\right],
\end{eqnarray}
where $V=\nabla \mu( r)/r$.
Denoting  $$\bar{V}=\sqrt{\varepsilon^2-\nu^2} \, V,~W=U+\nu  \, V,$$, we then have the following
\begin{eqnarray}
& \mathscr{E}^\nu(r,U)=  r \, \left(  \partial_t W  +  U \cdot \nabla W \right) +\nabla p(r) - \nu \,  \dive (\mu(r) \,  \nabla W)  \label{soldisQNS2} \\
 & -\dfrac{\nu}{2} \nabla (\lambda(r) \dive W) -  \sqrt{\eps^2-\nu^2}  \left(   \dive ( \mu(r) \,   \nabla \bar{V}) + \dfrac{1}{2} \nabla (\lambda(r) \dive \bar{V}) \right),  \nonumber  \\
& 0= r \, \left(  \partial_t \bar{V} + U  \cdot \nabla   \bar{V}  \right) +\sqrt{ \eps^2 - \nu^2} \,  \left(   \dive ( \mu(r)  \, ^t \nabla U)  + \dfrac{1}{2} \nabla(\lambda(r) \dive U \right).&  \label{soldisQNS3}
\end{eqnarray}
Before giving the definition let us recall that $\mathcal{E}_{{NSK}}(t)$ stands for
$$\mathcal{E}_{{NSK}}(t) =\mathcal{E}_{{NSK}}(\rho,\bar{v},w |  r, \bar{V},W)(t).$$
\begin{defi}\label{defsoldisI2}
Let us assume $\mu(\rho)=\rho^{(s+3)/2}$, $\gamma \geq s+2$ and $s \geq -1$.
Let $\rho_0$ and $u_0$ smooth enough.  The pair $(\rho,u) $ is a dissipative solution of \eqref{NSI0}-\eqref{NSI1}, \eqref{condI}
if the triplet $(\rho,\bar{v},w)$ (with ${\rho \, v=\nabla \mu( \rho)}, \, \bar{v}=\sqrt{\varepsilon^2-\nu^2} \, v,  \, w=u+\nu \, v$) satisfies 
\begin{eqnarray*}
\mathcal{E}_{{NSK}}(t) 
  \leq   \mathcal{E}_{{NSK}}(0) \, \exp(F^\nu t)+F^\nu  \int_0^t b_{{NSK}}(\xi) \, \exp(F^\nu \, (t-\xi)) \, d\xi
+b_{{NSK}}(t),\nonumber
\end{eqnarray*}
with $F^\nu$ given in {Corollary \ref{corsoldisNSQ}} and 
$$b_{{NSK}}(t)=\int_0^t \int_{\Omega} \left[\dfrac{\rho}{r} \mathscr{E}^\nu \cdot (W-w)\right]
$$
with $(r,\bar V, W)$ and $\mathscr{E}^\nu$ are defined as mentioned above from 
all given smooth function $U$.
\end{defi}

Noticing that each global weak solutions of the Navier-Stokes-Korteweg is global weak 
solutions of the augmented Navier-Stokes-Korteweg system, a direct consequence of the method is the following weak-strong uniqueness result.
 \begin{theorem}
 Let us assume $\mu(\rho)=\rho^{(s+3)/2}$, $\gamma \geq s+2$ and $s \geq -1$.
 Let us consider  $(\rho,u)$ a global weak solution to the compressible Navier-Stokes-Korteweg system and define 
 $w=u+\nu \,\nabla\mu(\rho)/\rho$ and 
 $\bar{v}= \sqrt{\eps^2-\nu^2} \, \nabla \mu(\rho)/\rho$.
    Let us assume that there exists $(r,U)$ a strong solution of  the compressible Navier-Stokes-Korteweg System and let us define
$W=U+\nu \nabla \mu(r)/r$ and $\bar{V}=\sqrt{\eps^2-\nu^2} \,{\nabla \mu(r)/r}.$  Assume that $(r,W, \bar{V})$ satisfies hypothesis \eqref{hypothesis}.
If $(\rho_0, u_0)=(r,U)(t=0)$  then $(\rho,\bar{v},w)=(r,\bar{V},W)$ or $(\rho,u)=(r,U)$, which corresponds to a  weak-strong uniqueness property.
\end{theorem}

Finally, let us give the definition \ref{defsoldisI2} in the particular case of $K(\rho)=1/\rho$ which corresponds to the quantum Navier-Stokes system. This one will be used in Section \ref{sec_limite}. To this end we introduce the function $\mathscr{E}_{NSQ}^\nu$ given by 
\begin{eqnarray}
\label{soldisQNS1}
\mathscr{E}_{NSQ}^\nu(r,U)= r \left( \partial_t U+ U \cdot  \nabla U  \right)+\nabla p(r)-2 \, \nu \,  \hbox{div}(rD(U))+\eps^2( \dive(r {}^t \nabla V ),
\end{eqnarray}
with $U$ a given smooth enough function, $r$ a strong solution of the mass equation \eqref{NSQF1} and $r V=\nabla \mu( r)$.
Defining  $$\bar{V}=\sqrt{\varepsilon^2-\nu^2} \, V,
\qquad W=U+\nu  \, V,$$ and Equation \eqref{NSQF1}
we obtain 
\begin{eqnarray}
&\mathscr{E}^\nu(r,U)=  r \, \left(  \partial_t W  +  U \cdot \nabla W \right) +\nabla p(r) - \nu \,  \dive (r \,  \nabla W) -  \sqrt{\eps^2-\nu^2} \dive ( r \,   \nabla \bar{V}) , \label{soldisQNS2} \\
& 0 = r \, \left(  \partial_t \bar{V} + U  \cdot \nabla   \bar{V}  \right) +\sqrt{ \eps^2 - \nu^2} \,  \dive ( r  \, ^t \nabla U).&  \label{soldisQNS3}
\end{eqnarray}
We define $\mathcal{E}_{{NSQ}}(t)$ by
$$\mathcal{E}_{{NSQ}}(t) =\mathcal{E}_{{NSK}}(t)~\mbox{with}~K(\rho)=1/\rho.$$
\begin{defi}\label{defsoldisI2bis}
Let $\rho_0$ and $u_0$ smooth enough.  The pair $(\rho,u) $ is a dissipative solution of \eqref{NSQI0}, \eqref{NSQI1}, \eqref{condI}
if the triplet $(\rho,\bar{v},w)$ (with ${\rho \, v=\nabla \mu( \rho)}, \, \bar{v}=\sqrt{\varepsilon^2-\nu^2} \, v,  \, w=u+\nu \, v$) satisfies 
\begin{eqnarray*}
\mathcal{E}_{{NSQ}}(t) 
  \leq   \mathcal{E}_{{NSQ}}(0) \, \exp(F^\nu t)+F^\nu  \int_0^t b_{{NSQ}}(\xi) \, \exp(F^\nu \, (t-\xi)) \, d\xi
+b_{{NSQ}}(t),\nonumber
\end{eqnarray*}
with $F^\nu$ given in {Corollary \ref{corsoldisNSQ}} and 
$$b_{{NSQ}}(t)=\int_0^t \int_{\Omega} \left[\dfrac{\rho}{r} \mathscr{E}_{NSQ}^\nu \cdot (W-w)\right]$$
with $(r, \bar V, W)$ and $\mathscr{E}_{NSQ}^\nu$ are defined as mentioned above from all given smooth function
$U$.
\end{defi}

\begin{rem}\label{remsoldis}
Note that by definition, using Corollary \ref{corsoldisNSQ},  all weak solution of \eqref{NSQI0}-\eqref{NSQI1}, \eqref{condI} is also a dissipative solution in the sense of Definition \ref{defsoldisI2bis}.
\end{rem}


 
\section{From the quantum Navier-Stokes system to the quantum Euler system: the viscous limit}\label{sec_limite}

We can now perform the limit of a dissipative solution of the quantum Navier-Stokes system to one of the quantum Euler system when the viscosity constant $\nu$ tends to zero. 
Thanks to the entropies, we have the following regularities on the global weak solution of the quantum Navier-Stokes equations:
$$\sqrt{\rho^\nu} \, \bar{v}^\nu  \in L^{\infty}(0,T;L^2(\Omega)),  
    \, \sqrt{\rho^\nu} \, w^\nu \in L^{\infty}(0,T;L^2(\Omega)), \,   H(\rho^\nu) \in L^{\infty}
(0,T,L^1(\Omega)),$$
where $$\bar{v}^\nu= \sqrt{\eps^2-\nu^2} \,  \nabla \log \rho^\nu, \quad w^\nu=u^\nu+ \nu \,   \nabla \log \rho^\nu.$$
The goal of this section is then to prove the following result:
\begin{theorem}
Let $\rho_0$ and $u_0$ smooth enough. Let $(\rho^\nu, u^\nu)$ be  a global weak solution to the quantum Navier-Stokes system \eqref{NSQI0}-{\eqref{NSQI1}} with initial conditions \eqref{condI}.
Let $(\rho, u)$ be the  weak limit  of $(\rho^{\nu}, u^{\nu})$  when $\nu$ tends to $0$ in the sense
$$\rho^\nu \rightharpoonup \rho \mbox{ weakly}\star \mbox{ in } L^\infty(0,T;L^\gamma(\Omega)) ,$$
$$\sqrt{\rho ^\nu} \, w^\nu \rightharpoonup  \sqrt \rho u \mbox{ weakly}\star \mbox{ in } L^{\infty}(0,T;L^2(\Omega)),$$
$$\sqrt{\rho^\nu} \, \bar{v}^\nu  \rightharpoonup \eps \,  \sqrt{ \rho} v
 \mbox{ weakly}\star \mbox{ in } L^{\infty}(0,T;L^2(\Omega)),$$
with  $\rho \, v=\nabla \rho$.
Then $(\rho, u)$  is a dissipative solution of the {quantum} Euler system \eqref{EQ1}-{\eqref{EQ2}} with initial conditions \eqref{condI}. 
\end{theorem}
\pr
According to Remark \ref{remsoldis}, the pair $(\rho^\nu,u^\nu)$ being an entropic weak solution, it is also a dissipative one. We want to prove that $(\rho,u)$, which is the limit of $(\rho^\nu,u^\nu)$ when $\nu$ tends to zero, is a dissipative solution of \eqref{EQ1}-\eqref{EQ2} satisfying the initial conditions  \eqref{condI}. The goal is then to prove that $(\rho,u)$ satisfies Definition \ref{Defi10}. 
  Let us define $v=\nabla \log \rho$ (because in this case $\mu(\rho)=\rho$).
  Let $U$ be smooth function and  let $(r,\mathscr{E}(r,U),V)$ be defined with $V=\nabla \log r$, \eqref{NSQF1} and  \eqref{soldisI2}.
   We define
$$\bar{V}^\nu=\sqrt{\varepsilon^2-\nu^2} \, V, \quad ~W^\nu=U+\nu \, V.$$
Then it is easy to see that $(r,\bar{V}^\nu, W^\nu)$ is candidate for \eqref{NSQF1}, \eqref{soldisQNS2}-\eqref{soldisQNS3} with
$$\mathscr{E}_{NSQ}^\nu(r,U)=\mathscr{E}(r,U)-2 \, \nu \,  \hbox{div}(rD(U)).$$
Then using Definition \ref{defsoldisI2bis}, $(\rho^\nu,u^\nu)$ being a dissipative solution we have
\begin{equation} \label{proofconv1}
\mathcal{E}_{NSQ}(t) 
  \leq   \mathcal{E}_{NSQ}(0) \, \exp(F^\nu t)+F^\nu  \int_0^t b_{NSQ}^\nu(\xi) \, \exp(F^\nu \, (t-\xi)) \, d\xi
+b_{NSQ}^\nu(t),
\end{equation}
with
\begin{eqnarray*}
F^\nu & =   & C\left(1 +   \dfrac{\nu}{\eps^2-\nu^2 }\right),
\end{eqnarray*}
and
$$b_{NSQ}^\nu(t)=\int_0^t \int_{\Omega} \left[\dfrac{\rho}{r} \left(\mathscr{E} -2 \, \nu \,  \hbox{div}(rD(U))\right)\cdot (W-w)\right].$$
Since by definition we have
\begin{eqnarray}
&{}&\mathcal{E}_{NSQ}(\rho^\nu,\bar{v}^\nu,w^\nu |  r, \bar{V}^\nu,W^\nu)(t)=
\dfrac{1}{2} \int_{\Omega}
\rho^\nu \, \left(  |\bar{v}^\nu-\bar{V}^\nu|^2 +  |w^\nu-W^\nu|^2  \right) \nonumber\\
 & +& \int_{\Omega} \left( H(\rho^\nu)-H(r)-H'(r)(\rho^\nu-r) \right)
+\nu \, \int_0^t \,  \int_{\Omega} \rho^\nu \,  \left(|\nabla \bar{v}^\nu- \nabla \bar{V}|^2 +  |\nabla w^\nu-\nabla W |^2\right), \nonumber
\end{eqnarray}
we easily obtain
\begin{eqnarray*}
\dfrac{1}{2} \int_{\Omega}
\rho^\nu \, \left(  |\bar{v}^\nu-\bar{V}^\nu|^2 +  |w^\nu-W^\nu|^2  \right) 
  &+ &\int_{\Omega} \left( H(\rho^\nu)-H(r)-H'(r)(\rho^\nu-r) \right) \\
&\leq& \mathcal{E}_{NSQ}(\rho^\nu,\bar{v}^\nu,w^\nu |  r, \bar{V}^\nu,W^\nu)(t)  
\end{eqnarray*}
and
\begin{eqnarray*}
\mathcal{E}_{NSQ}(\rho^\nu,\bar{v}^\nu,w^\nu |  r, \bar{V}^\nu,W^\nu)(0)=
\dfrac{1}{2} \int_{\Omega}
\rho^\nu \, \left(  |\bar{v}^\nu-\bar{V}^\nu|^2 +  |w^\nu-W^\nu|^2  \right) (0) \\
  + \int_{\Omega} \left( H(\rho^\nu)-H(r)-H'(r)(\rho^\nu-r) \right)(0).
\end{eqnarray*}
Then \eqref{proofconv1} gives
\begin{eqnarray}
&{}&\dfrac{1}{2} \int_{\Omega}
\rho^\nu \, \left(  |\bar{v}^\nu-\bar{V}^\nu|^2 +  |w^\nu-W^\nu|^2  \right) (t)
  +\int_{\Omega} \left( H(\rho^\nu)-H(r)-H'(r)(\rho^\nu-r) \right)(t)  \nonumber\\
  & \leq &   \mathcal{E}_{NSQ}(\rho^\nu,\bar{v}^\nu,w^\nu |  r, \bar{V}^\nu,W^\nu)(0)
 \exp  (F ^\nu\, t)  +F^\nu  \int_0^t b^\nu_{NSQ}(\xi) \, \exp(F^\nu \, (t-\xi)) \, d\xi +b^\nu_{NSQ}(t). \nonumber
\end{eqnarray}
  It remains now to pass to the limit $\nu$ tends to zero in this inequality. 
Clearly, using the lower semi-continuity of the term $\mathcal{E}_{NSQ}(\rho^\nu,\bar{v}^\nu,w^\nu |  r, \bar{V}^\nu,W^\nu)$, the left-hand side is greater than
\begin{eqnarray*}
\dfrac{1}{2} \int_{\Omega}
\rho \, \left(  \varepsilon^2 | \, v-V|^2 +  |u-U|^2  \right) (t) + \int_{\Omega}
H(\rho | r)(t), 
\end{eqnarray*}
which is $\mathcal{E}_{EuQ}(\rho,u,v|r,U,V)(t)$ ({\it i.e.} $\mathcal{E}_{EuK}(\rho,u,v|r,U,V)(t)$ given by \eqref{relentEK} with $K(\rho)=1/\rho$). 
  For the right hand side, we use the direct limit of the term 
  $\mathcal{E}_{NSQ}(\rho^\nu,\bar{v}^\nu,w^\nu |  r, \bar{V}^\nu,W^\nu)(0)$ (through the expression of the initial data) and
 $b^\nu_{NSQ}$ tends to
$$b_{EuQ}(t)=\int_0^t \int_{\Omega} \left[\dfrac{\rho}{r} \mathscr{E} \cdot (U-u)\right],$$
to conclude that
\begin{eqnarray*}
\mathcal{E}_{EuQ}(t)    \leq    \mathcal{E}_{EuQ}(0)
 \exp  (C \, t)
 + b_{EuQ}(t) +C   \, \displaystyle  \int_0^t \exp(C \, (t-\xi))b_{EuQ}(\xi) \, d\xi,
\end{eqnarray*}
where $C=C(\eps^2,r,U,V)$ is a uniformly bounded constant on $\R^+ \times \Omega$.
Therefore we finally obtain that $(\rho,u)$ satisfies the Definition \ref{Defi10} and then is a dissipative solution of \eqref{EQ1}-\eqref{EQ2}, \eqref{condI}.
 \cqfd



\section{Appendix}
\vskip0.5cm



\subsection{Technical lemmas on modulated quantities}\label{sec_lemtech}
In this section we give some technical lemmas  which are used in the paper. 

\medskip

\noindent We  introduce the function $\phi$ defined by
 \begin{equation} \label{phi}
 \phi(\tau)= \int_0^\tau \dfrac{p'(\mu^{-1}(s))}{\mu'(\mu^{-1}(s))} ds,
 \end{equation}
and the two functions
\begin{equation}\label{phi1}
\phi_1(\rho|r)=\phi(\mu(\rho))-\phi(\mu(r))-\phi'(\mu(r))(\mu(\rho)-\mu(r)),
\end{equation}
\begin{equation}\label{phi2}
\phi_2(\rho|r)=\phi''(\mu(r))(\mu(\rho)-\mu(r)) \, r  - \rho \, (\phi'(\mu(\rho))-\phi'(\mu(r))).
\end{equation}

\medskip

\begin{rem}
Note that in the case $K(\rho)=1/\rho$, which gives {\rm (}using \eqref{muK}{\rm )} $\mu(\rho)=\rho$, these two functions are directly linked to $H(\rho|r)$. Indeed, in this case we have
\begin{eqnarray*}
\phi_1(\rho|r) &  = &  p(\rho)-p(r)-p'(r)(\rho-r)=(\gamma-1)H(\rho|r),\\
\phi_2(\rho|r) & =  & \rho p'(r)-\rho p'(\rho)+r p''(r)(\rho-r) =-\gamma (\gamma-1)H(\rho | r).
\, \square
\end{eqnarray*}
\end{rem}

\medskip

As usually in compressible flows (see \cite{FeNo}) let us define the set $\mathcal{F}$ by
$$\mathcal{F}=\left\{ \rho \leq \dfrac{r}{2} ~\hbox{or}~\rho \geq 2r\right\}.$$

 Let us now give some technical lemmas which will be used in the following. First of all, following \cite{FNL} we have
\begin{lem}\label{lem_tech1}
Assuming $p$ smooth, $p(0)=0$, $p'(\rho)>0  \quad \forall \rho >0$,  $\displaystyle \lim_{\rho \rightarrow \infty}  \dfrac{p'(\rho)}{\rho^{\alpha-1}} =a >0$ for  $\alpha  >1$,
 we have:
 \begin{eqnarray*} 
H(\rho | r) \geq C(r) (\rho-r)^2 \mbox{ if } \rho \in \mathcal{F}^c~\hbox{ and }~
H(\rho | r) \geq C(r) (1+\rho)^\gamma \hbox{ otherwise},
  \end{eqnarray*}
 with $C(r)$ uniformly bounded for $r$ belonging to compact sets in $\R^+ \times \Omega$.
\end{lem}

Concerning the functions $\phi_1$ and $\phi_2$, we can show
\begin{lem}\label{lem_tech2}
Let us assume that $\mu(\rho)=\rho^{(s+3)/2}$ with $\gamma \geq s+2$ and $s \geq -1$. Assume $\phi_i$ with $i=1,2$ defined by \eqref{phi}--\eqref{phi2}. Then 
\begin{eqnarray*}
|\phi_i(\rho|r)| \leq C(r) |\rho-r|^2  \mbox{ if } \rho \in \mathcal{F}^c~\hbox{ and }~
|\phi_i(\rho|r)| \leq C(r)(1+\rho)^\gamma \hbox{ otherwise},
\end{eqnarray*}
with $C(r)$ uniformly bounded for $r$ belonging to compact sets in $\R^+ \times \Omega$.
\end{lem}

\begin{rem}
Let us remark that the choice $\mu(\rho)=\rho^{(s+3)/2}$  with $s \in \R$ and the assumption $\gamma \geq 2+s$ correspond to the case considered in \cite{GLT} because $K(\rho)$ is of order $\rho^s$. Moreover, for the particular case of interest in this paper $K(\rho)=1/\rho$ (\it{i.e.} $s=-1$), the assumption $2+s \leq \gamma$ is trivially satisfied since we have $\gamma > 1$.
\end{rem}

\noindent{\it Proof of the lemma for $\phi_1$.}
Using Taylor expansions and the fact that $\phi''(\mu(c))$, $\mu'(c)$ are bounded with  $c$ in a compact we easily obtain
\begin{eqnarray*}
|\phi_1(\rho|r)| &  \leq  &  C(r) | \mu(\rho)-\mu(r)|^2  \leq  C(r) |\rho-r|^2 ~\hbox{on}~\mathcal{F}^c.
 \end{eqnarray*}
Moreover, since $$p(\rho)=\rho^\gamma, \, K(\rho)=\dfrac{(s+3)^2}{4} \rho^s, \, \mu'(\rho)=\sqrt{\rho K(\rho)},$$ we have $\phi(\tau)=\tau^{2\gamma/(s+3)}$ and then by definition
\begin{eqnarray*}
|\phi_1(\rho|r)| & = & \left| \rho^\gamma -r^{\gamma} - \dfrac{2\gamma}{s+3}r^{\frac{2 \gamma-(s+3)}{2}} (\rho ^\frac{s+3}{2}-r^\frac{s+3}{2})\right|,
\end{eqnarray*}
which gives
\begin{eqnarray*}
|\phi_1(\rho|r)|  \leq  C(r)(1+ \rho)^{\gamma}~\hbox{on}~\mathcal{F},
\end{eqnarray*}
since by assumption $2\gamma \geq 2(s+2) \geq s+3$ with $s \geq -1$.

\noindent {\it Proof of the lemma for $\phi_2$.}
Let us write  $\theta=\dfrac{s+3}{2}$ then $\mu(\rho)= \rho^{\theta}$ and $\phi(\rho) = \rho^{\gamma/\theta}$. Then 
\begin{eqnarray*}
\phi_2(\rho|r) & =  & \dfrac{2 \gamma}{s+3} \left[ \left( \dfrac{2 \gamma}{s+3}-1 \right) r^{\gamma-s-3}(\rho^{\theta}-r^{\theta}) r - \rho(\rho^{\gamma-\theta}-r^{\gamma-\theta}) \right] \\
 & = &  \dfrac{2 \gamma}{s+3} \left[ \left( \dfrac{2 \gamma}{s+3}-1 \right) r^{\gamma-s-2}(\rho^{\theta}-r^{\theta}) - \rho^{1+\gamma-\theta}+\rho \, r^{\gamma-\theta} \right] \\
 & = &  \dfrac{2 \gamma}{s+3} \left[ \left( \dfrac{2 \gamma}{s+3}-1 \right) r^{\gamma-s-2}(\rho^{\theta}-r^{\theta})  -f(\mu(\rho))+\rho \, r^{\gamma-\theta} \right]
\end{eqnarray*}
with $f(\rho)= \rho ^{\frac{\gamma+1}{\theta}-1} $. Note that we have
 \begin{eqnarray*}
f(\rho|r) & = &  f(\mu(\rho))-f(\mu(r))-f'(\mu(r))(\mu(\rho)-\mu(r)) \\
& =  & f(\mu(\rho))-(r^{\theta})^{\frac{\gamma+1}{\theta-1}} -\left( \dfrac{\gamma+1}{\theta}-1 \right)  (r^{\theta})^{\frac{\gamma+1}{\theta}-2}(\rho^{\theta}-r^{\theta}). \\
\end{eqnarray*}
Then
\begin{eqnarray*}
\phi_2(\rho|r) & =  & \dfrac{2 \gamma}{s+3} \left[ \left (\dfrac{2 \gamma}{s+3}-1 \right) r^{\gamma-s-2}(\rho^{\theta}-r^{\theta})  -f(\rho|r) -  (r^{\theta})^{\frac{\gamma+1}{\theta}-1} \right] \\
&   &  - \dfrac{2 \gamma}{s+3} \left[ \left( \dfrac{2(\gamma+1)}{s+3}-1 \right)  (r^{\theta})^{\frac{\gamma+1}{\theta}-2}(\rho^{\theta}-r^{\theta}) - \rho \, r ^{\gamma -\theta}  \right] \\
& =  & 
\dfrac{2 \gamma}{s+3} \left[ \left( \dfrac{2 \gamma}{s+3}-1 \right) r^{\gamma-s-2}(\rho^{\theta}-r^{\theta})  -f(\rho|r) - r^{1+\gamma-\theta} \right] \\
& &  - \dfrac{2 \gamma}{s+3} \left[\left( \dfrac{2(\gamma+1)}{s+3}-1 \right)  (r^{\gamma-s-2}(\rho^{\theta}-r^{\theta}) -  \rho \, r ^{\gamma -\theta}  \right] \\
 & = & \dfrac{2 \gamma}{s+3} \left[ -\dfrac{1}{\theta}  r^{\gamma-s-2}(\rho^{\theta}-r^{\theta})  -f(\rho|r)  + \rho \, r ^{\gamma -\theta} - r^{1+\gamma-\theta}  \right]. 
\end{eqnarray*}
This can be written $\phi_2(\rho|r)=\dfrac{2 \gamma}{s+3} (-f(\rho|r)+g(\rho|r))$ with
\begin{eqnarray*}
g(\rho|r) & = & r^{\gamma-\theta}\left[ \rho-r-\dfrac{1}{\theta} r^{\tau-s-2}(\rho^{\theta}-r^{\theta}) \right] \\
 & = &  r^{\gamma-\theta} \left[ ( \rho^{\theta})^{1/\theta}-(r^{\theta})^{1/\theta}-\dfrac{1}{\theta}(r^{\theta})^{1/\theta-1}( \rho^{\theta}-r^{\theta}) \right].
 \end{eqnarray*}
In the case  $\rho \in \mathcal{F}^c$, using Taylor expansions this leads to
$$|f(\rho|r)| \leq C(r) |\mu(\rho)-\mu(r)]^2 \leq C(r) |\rho-r|^2,$$
$$|g(\rho|r)| \leq  C(r)  |\rho^\theta-r^\theta|^2 \leq  C(r)  |\rho-r|^2,$$
and then
$$|\phi_2(\rho|r)| \leq  C(r) |\rho-r|^2.$$
When $\rho \in \mathcal{F}$, since $2 \gamma \geq 2s+4 \geq s+3$ and $s+3 \geq 2$,
\begin{eqnarray*}
|\phi_2(\rho|r)|  & \leq & C(r)| r^{\gamma-(s+3)}(\rho^{\frac{s+3}{2}}-r^{\frac{s+3}{2}})- \rho  (\rho^{\gamma-\frac{s+3}{2}}-r^{\gamma-\frac{s+3}{2}})|  \leq  C(r) (1+\rho)^{\gamma}.
\end{eqnarray*}
This completes the proof of Lemma \ref{lem_tech2}.
\cqfd

Using Lemmas \ref{lem_tech1} and \ref{lem_tech2}, we directly obtain
\begin{lem}\label{lem_tech3}
Let us assume that $\mu(\rho)= \rho^{(s+3)/2}$ with $\gamma \geq s+2$ and $s \geq -1$. We have
\begin{eqnarray*}
|\phi_1(\rho|r)| \leq C(r) H(\rho|r) \quad \hbox{and} \quad |\phi_2(\rho|r)| \leq C(r) H(\rho|r),
\end{eqnarray*}
with $C(r)$ uniformly bounded for $r$ belonging to compact sets in $\R^+ \times \Omega$.
\end{lem}

Let us now prove the following lemma
\begin{lem}\label{lem_tech4}
Let us assume that $\mu(\rho)= \rho^{(s+3)/2}$ with $\gamma \geq s+2$ and $s \geq -1$. We have
\begin{eqnarray*}
\rho|\mu'(\rho)-\mu'(r)|^2 \leq C(r) H(\rho|r),
\end{eqnarray*}
with $C(r)$ uniformly bounded for $r$ belonging to compact sets in $\R^+ \times \Omega$.
\end{lem}

\proof
$$\rho | \mu'(\rho)-\mu'(r) |^2 =  \rho | \mu'(\rho)-\mu'(r) |^2 1_{\mathcal F} + \rho | \mu'(\rho)-\mu'(r) |^2 1_{\mathcal F^c}. $$
We have
\begin{eqnarray*}
 \rho | \mu'(\rho)-\mu'(r) |^2 1_{\mathcal F}  \leq  2   \rho (| \mu'(\rho)|^2+|\mu'(r) |^2 ) 1_{\mathcal F}  \leq  \dfrac{(s+3)^2}{2}  \rho ^{s+2} 1_{\mathcal F}+ 2 C(r) \rho \,   1_{\mathcal F} .
\end{eqnarray*}
Using $ \rho ^{s+2} \leq (1+ \rho )^{s+2}$ and the assumption $\gamma \geq s+2$ in the first term, and, the assumption $\gamma >1$ in the second one, we obtain:
\begin{eqnarray*}
\rho | \mu'(\rho)-\mu'(r) |^2 1_{\mathcal F} \leq \dfrac{(s+3)^2}{2}  (1+\rho)^\gamma 1_{\mathcal F}+ 2C(r) (1+\rho)^\gamma 1_{\mathcal F}  \leq  C(r) (1+\rho)^\gamma 1_{\mathcal F}.
 \end{eqnarray*}
Moreover,
 \begin{eqnarray*}
 \rho | \mu'(\rho)-\mu'(r) |^2 1_{\mathcal F^c} = \dfrac{s+3}{2} \rho \left| \rho ^{\frac{s+1}{2}}-r^{\frac{s+1}{2}}\right|^2 1_{\mathcal F^c} = \dfrac{s+3}{2}  \rho \frac{ | \rho ^{s+1}-r^{s+1} |^2} { | \rho ^{\frac{s+1}{2}}+r^{\frac{s+1}{2}} |^2}1_{\mathcal F^c} , 
\end{eqnarray*}
and then
\begin{eqnarray*}
 \rho | \mu'(\rho)-\mu'(r) |^2 1_{\mathcal F^c}  \leq \dfrac{s+3}{2} \rho \frac{ | \rho ^{s+1}-r^{s+1} |^2} { | r ^{\frac{s+1}{2}} |^2}1_{\mathcal F^c}  \leq C(r)  | \rho ^{s+1}-r^{s+1} |^2 1_{\mathcal F^c}   \leq C(r)  | \rho-r  |^2 1_{\mathcal F^c}. \\
 \end{eqnarray*}
Using lemma \ref{lem_tech1}, we finally obtain the result.
\cqfd

\medskip

\noindent {\it An important relation.}
The last technical and important lemma is
\begin{lem}\label{lem_tech5}
Let us assume that $\mu(\rho)= \rho^{(s+3)/2}$ with $\gamma \geq s+2$ and $s \geq -1$.
 We have
 \begin{eqnarray*}
\rho\left(\dfrac{p'(\rho)}{\mu'(\rho)} v -\dfrac{p'(r)}{\mu'(r)}V \right)\cdot(v-V)=
[\nabla \phi_1(\rho|r) +\phi_2(\rho|r)V]\cdot V + \rho \frac{p'(\rho)}{\mu'(\rho)} |V-v|^2
\end{eqnarray*}
with $\phi_1$ and $\phi_2$ defined by  \eqref{phi}--\eqref{phi2}.
\end{lem}

\begin{rem} This lemma generalizes to general $\mu(\rho)$  the relation (5) established in \cite{BNV} when $\mu(\rho)= \rho$.  This is an important lemma which  helps to control the terms coming from the pressure in the relative entropy at the Navier-Stokes level.
\end{rem}

\proof 
Remark first that 
$$\rho\left(\dfrac{p'(\rho)}{\mu'(\rho)} v -\dfrac{p'(r)}{\mu'(r)}V \right)\cdot(v-V)
=  \rho \frac{p'(\rho)}{\mu'(\rho)} |V-v|^2 
 + \rho \, \left( \dfrac{p'(\rho)}{\mu'(\rho)}-\dfrac{p'(r)}{\mu'(r)} \right) (v-V)\cdot V
$$

We have
\begin{eqnarray*}
&{}& \rho \, \left( \dfrac{p'(\rho)}{\mu'(\rho)}-\dfrac{p'(r)}{\mu'(r)} \right) (v-V) \\
  & = & \left( \dfrac{p'(\rho)}{\mu'(\rho)}-\dfrac{p'(r)}{\mu'(r)} \right) (\rho \, v- \rho \, V)  \\\
   & = & \left( \dfrac{p'(\rho)}{\mu'(\rho)}-\dfrac{p'(r)}{\mu'(r) } \right) \nabla (\mu(\rho)) -  \dfrac{\rho}{r} \left( \dfrac{p'(\rho)}{\mu'(\rho)}-\dfrac{p'(r)}{\mu'(r)}  \right) \nabla (\mu(r)).
\end{eqnarray*}
Moreover, it is easy to see that by definition
\begin{eqnarray*}
\nabla (\phi_1(\rho|r)) & = &  \phi'(\mu(\rho)) \nabla \mu(\rho) -\phi''(\mu(r))(\mu(\rho)-\mu(r))\nabla\mu(r) - \phi'(\mu(r))  \nabla\mu(\rho) \\
& = & \left(\dfrac{p'(\rho)}{\mu'(\rho)}-\dfrac{p'(r)}{\mu'(r)}\right)\nabla \mu(\rho)-\phi''(\mu(r))(\mu(\rho)-\mu(r))\nabla\mu(r),
\end{eqnarray*}
and then using the definition of $\phi_2(\rho|r)$,
\begin{eqnarray*}
\rho\left(\dfrac{p'(\rho)}{\mu'(\rho)}-\dfrac{p'(r)}{\mu'(r)}\right)(v-V)=\nabla \phi_1(\rho|r) +\phi_2(\rho|r)V.
\end{eqnarray*}
\cqfd

\subsection{Equivalence of $\mathcal{E}_{EuK}$ and the relative entropy in \cite{GLT}}{\label{annex_equivalence}}

Let us consider the relative entropy functional, denoted $\mathcal{E}_{EuK}(\rho,u,v | r,U,V)$ and defined by \eqref{relentEK}. The goal of this section is to prove that this relative entropy 
is equivalent  to the relative entropy defined by (2.23) in \cite{GLT} under  the concavity assumption on $K$ with $K(\rho)= \rho^s$. Let us first recall the relative entropy $\mathcal{E}_{EuK}^{GLT}$ defined in \cite{GLT}. It reads
\begin{eqnarray}
\label{relentEKGLT}
&{}&\mathcal{E}_{EuK}^{GLT}(\rho,u,\nabla \rho \vert r,U, \nabla r)=
  \frac{1}{2} \int_\Omega \rho |u-U|^2+  \dfrac{1}{2} \varepsilon^2 \int_\Omega I_T +\int_\Omega H(\rho | r),  
\end{eqnarray}
where $$I_T= K(\rho) |\nabla \rho|^2 -K(r) |\nabla r|^2 -K'(r) |\nabla r|^2 (\rho-r) -2 K(r) \nabla r (\nabla \rho-\nabla r).$$
Note that $I_T$ corresponds to the term $K(\rho)|\nabla\rho|^2$ linearized in the variables
$(\rho,q)$ where $q=\nabla\rho$.
 Let us now introduce the quantity
\begin{eqnarray*}
I_{EuK}^2  & = &  \rho \left|\sqrt{\frac{K(\rho)}{\rho}} \nabla \rho - \sqrt{\frac{K(r)}{r}} \nabla r\right|^2  =  \left|\sqrt{K(\rho)} \nabla \rho- \sqrt{\frac{\rho}{r}} \sqrt{K(r)} \nabla r\right|^2.
  \end{eqnarray*}
  Then our Euler-Korteweg modulated energy reads
 $$\mathcal{E}_{EuK}(\rho,u,v | r,U,V)= \dfrac{1}{2} \int \rho|u-U|^2 +\dfrac{1}{2}  \eps^2 \int  I^2_{EuK} + \int H(\rho|r),$$
 where $v= \sqrt{K(\rho)}\nabla \rho/\sqrt \rho$ and $V= \sqrt{K(r)}\nabla r/\sqrt r$.
Let us prove that  under the hypothesis on $K$ introduced in \cite{GLT}
$$\mathcal{E}_{EuK}(\rho,u,v | r,U,V)=0 \qquad  \Leftrightarrow  \qquad 
    \mathcal{E}_{EuK}^{GLT}(\rho,u,\nabla \rho \vert r,U, \nabla r) = 0.$$
  If so, we prove by this way that our relative entropy and the one in \cite{GLT}
 are equivalent under the hypothesis in \cite{GLT}. Our convergence result will therefore be more general that the one in \cite{GLT} because it does not asked for concavity hypothesis on $K(\rho)$. First let us prove the following lemma:

\medskip
\begin{lem}  \label{identity} 
We have the equality
$$I_{\rm EuK}^2+ K(r) \,  |\nabla r|^2 \, \left| \sqrt{\frac{K(r)}{K(\rho)}}- \sqrt{\frac{\rho}{r}}\right|^2- K(r)^2 \,  |\nabla r|^2  \, \left( \frac{1}{K(\rho)}-\frac{1}{K(r)}+\frac{K'(r)}{K(r)^2} (\rho-r)\right) $$
$$ \hskip3cm  = I_T + 2 \sqrt{K(r)} \nabla r I_{\rm EuK} \
    \left( \sqrt{\frac{K(r)}{K(\rho)}}- \sqrt{\frac{\rho}{r}}\right)
$$
\end{lem}

\noindent {\bf Proof.}
After computations, we check that
\begin{eqnarray*}
I_{EuK}^2 & = & I_T + K'(r) |\nabla r|^2 (\rho-r) +\frac{\rho}{r} K(r) |\nabla r|^2 \\
 & &  - 2 \sqrt{\frac{\rho}{r} } \sqrt{K(\rho)} \sqrt{K(r)} \nabla \rho \cdot \nabla r 
+2K(r) \nabla r \cdot \nabla \rho-K(r) |\nabla r|^2 \\
 & =  & I_T+I_1
 \end{eqnarray*}
where
$$I_1=K(r)|\nabla r|^2\left(\dfrac{\rho}{r}-1+\dfrac{K'(r)}{K(r)}(\rho-r)+2 \frac{\sqrt{K(r)}}{\sqrt{K(\rho)}} \sqrt{\frac{\rho}{r}}-2 \frac{\rho}{r} \right) +I_2,$$
with
$$I_2= 2 \,  \sqrt{K(r)} \, \nabla r  \, I_3 \,  I_{EuK} \qquad 
      \hbox{ and } \qquad  I_3= \sqrt{\frac{K(r)}{K(\rho)}}- \sqrt{\frac{\rho}{r}}.$$

\begin{cor} Let $K(\rho) = \rho^s$ with $-1\le s \le 0$, then 
$$\mathcal{E}_{EuK}^{GLT}(\rho,u,\nabla \rho \vert r,U, \nabla r) = 0 
\qquad  \Leftrightarrow  \qquad 
  \mathcal{E}_{EuK}(\rho,u,v | r,U,V)=0.$$
\end{cor}

\noindent {\bf Proof.} Under the assumption on $K$, we check that
 \begin{eqnarray*}
 I_3^2  & = & \left( \sqrt{\dfrac{K(r)}{K(\rho)}}-\sqrt{\dfrac{\rho}{r}} \right )^2 \\
 & = &\left( \sqrt{\left( \dfrac{\rho}{r} \right)^{-s}}-\sqrt{\dfrac{\rho}{r}} \right)^2 \\
 & \leq  & 2 \left( \sqrt{\left( \dfrac{\rho}{r} \right)^{-s}}-1 \right)^2+2 \left( 1-\sqrt{\dfrac{\rho}{r}} \right)^2 \\
 &  \leq  & 2 \dfrac{1}{r^{-s} } (\sqrt{\rho^{-s}}-\sqrt{r^{-s}})^2 +\dfrac{2}{r} (\sqrt{r}-\sqrt{\rho})^2 \\
  & \leq & 
 \dfrac{2}{r^{-s}} | \rho^{-s}-r^{-s}|+\dfrac{2}{r} |r-\rho|
 \end{eqnarray*}
 with $0\le -s \le 1$. Assume   $\mathcal{E}_{EuK}^{GLT}(\rho,u,\nabla \rho \vert r,U, \nabla r) = 0$, then $I_3=0$ and 
 $$\left( \frac{1}{K(\rho)}-\frac{1}{K(r)}+\frac{K'(r)}{K(r)^2} (\rho-r)\right)=0.$$ 
 Therefore using Lemma \ref{identity} we conclude 
 $ \mathcal{E}_{EuK}(\rho,u,v | r,U,V)=0$ (the inverse follows the same lines).
  This ends  the proof. \cqfd

\vskip0.5cm
\subsection{Definition of the operators}{\label{annex_operators}}
For the convenience of the reader we recall in this Section all the definitions of the operators used in this article. The definitions used here are the ones presented in \cite{BoFa06} in Appendix A.

\bigskip
Let $f$ be a scalar, $u, v$ two vectors and $\sigma=(\sigma_{ij})_{1 \leq i,j \leq d}$ a tensor field defined on $\Omega \subset \R^d$ smooth enough. 

\begin{itemize}
	\item Denoting by $v_1, \cdots, v_d$ the coordinates of $v$, we call {\it divergence} of $v$ the scalar given by:
	$$\hbox{div} (v)= \sum_{i=1}^d \dfrac{\partial v_i}{\partial x_i}.$$
	\item We call {\it laplacian} of $f$ the scalar given by:
	$$\Delta f=\hbox{div} (\nabla f)= \sum_{i=1}^d \dfrac{\partial^2f}{\partial x_i^2}.$$
	\item We call {\it gradient} of $v$ the tensor given by:
	$$\nabla v =\left(\dfrac{\partial v_i}{\partial x_j}\right)_{1\leq i,j \leq d}. $$
	\item We call {\it divergence} of $\sigma$ the vector given by:
	$$\hbox{div}(\sigma) =\left(\sum_{j=1}^d \dfrac{\partial \sigma_{ij}}{\partial x_j}\right)_{1\leq i \leq d} .$$
	\item We call {\it laplacian} of $v$ the vector given by:
	$$\Delta v =\hbox{div}(\nabla v).$$
	\item We call {\it tensor product} of $u$ and $v$ the tensor given by:
	$$u \otimes v = \left(u_i v_j\right)_{1\leq i,j \leq d}. $$
\end{itemize}

\begin{prop}\label{propAPA}
Let $u,v,w$ three smooth enough vectors on $\Omega$ and $r$ a scalar smooth enough on $\Omega$. We have the following properties.
\begin{itemize}
	\item $(u \otimes v)w=(v \cdot w)u$,
	\item $\dive (u \otimes v)=(\dive v)u+(v\cdot \nabla)u$,
	\item $\dive(r \, u)=\nabla r \cdot u+ r \, \dive u$,
	\item $\dive(r \, u\otimes v)=(\nabla r \cdot v)u + r (v \cdot \nabla)u+r \, \dive(v)u$.
\end{itemize}
\end{prop}

\begin{defi}
Let $\tau$ and $\sigma$ be two tensors of order 2. We call scalar product of the two tensors the real defined by:
$$\sigma : \tau = \sum_{1 \leq i,j \leq d} \sigma_{ij}\tau_{ij}.$$
The norm associated to this scalar product is simply denoted by $| \cdot |$ in such a way that 
$$|\sigma|^2=\sigma : \sigma.$$
\end{defi}

\begin{rem}\label{remappA}
By definition we have
$$ \sigma : \tau = {}^t\sigma : {}^t\tau$$
\end{rem}


\section*{Acknowledgements}
The third author acknowledges support from the team INRIA/RAPSODI and the Labex CEMPI (ANR-11-LABX-0007-01). The first author acknowledges the project  TELLUS INSU-INSMI "Approche crois\'ee pour fluides visco-\'elasto-plastiques: vers une meilleure compr\'ehension des zones solides/fluides". The first and the second authors acknowledge the ANR Project FRAISE managed by C. {\sc Ruyer-Quil}.

\bigskip

\noindent {\bf Conflict of interest:} The authors have no conflicts of interest to declare.



\end{document}